\documentclass[11pt]{article}
\usepackage{amssymb, amsmath, amsthm}
\usepackage{epsfig}
\usepackage{color}
\usepackage{graphics}
\numberwithin{equation}{section}
\newtheorem{teo}{Theorem}[section]
\newtheorem{cor}[teo]{Corollary}
\newtheorem{prop}[teo]{Proposition}
\newtheorem{lemma}[teo]{Lemma}

\newtheorem{remark}[teo]{Remark}
\newtheorem{example}[teo]{Example}
\newtheorem{definition}[teo]{Definition}

\theoremstyle{definition}

\theoremstyle{remark}

\font\svfilt=msbm7

\font\filt=msbm10

\def\neweq#1{\begin{equation}\label{#1}}
\def\endeq{\end{equation}}

\def\proof{\noindent{\it Proof.} }

\def\A{\mathcal A}
\def\QC{\mathcal Q ^n}

\def\sfe{{\mathbb S}^{n-1}}
\def \cl {\mathop {\rm cl}\nolimits}

\def\H{{\mathcal H}}
\def\K{\mathcal K}
\def\L{\mathcal L}

\def\Per{{\rm Per}}

\def \relint {\mathop {\rm int}\nolimits}
\def \R {\hbox{\filt R}}

\def \sR {\hbox{\svfilt R}}

\def\bee{\begin{eqnarray*}}
\def\ene{\end{eqnarray*}}
\def\ep{\varepsilon}

\def\be{\begin{equation}}
\def\en{\end{equation}}

\begin{document}

\title{\bf Quermassintegrals of quasi-concave functions \\ 
and generalized Pr\'ekopa-Leindler inequalities}

\author{S. G. BOBKOV\footnote{Research is partially supported by
NSF grant and Simons Fellowship}, \ A. COLESANTI, \ I. FRAGAL\`A}

\date{}
\baselineskip14pt \maketitle

\begin{abstract}
\hskip-5mm 
We extend to a functional setting the concept of quermassintegrals, well-known within the Minkowski theory of convex bodies.
We work in the class of quasi-concave functions defined on the Euclidean space, and
with the hierarchy of their subclasses given by $\alpha$-concave functions.
In this setting, we investigate the most relevant features of functional quermassintegrals,
and we show they inherit the basic properties of their classical geometric counterpart.
As a first main result, we prove a Steiner-type formula which holds true
by choosing a suitable functional equivalent of the unit ball.
Then, we establish concavity inequalities for quermassintegrals and for other general hyperbolic functionals,
which generalize the celebrated Pr\'ekopa-Leindler and Brascamp-Lieb inequalities.
Further issues that we transpose to this functional setting are:
integral-geometric formulae of Cauchy-Kubota type, valuation property and isoperimetric/Uryshon-like inequalities. 
\end{abstract}

\medskip
\noindent {\small {\sl 2010\,MSC}: 28B, 46G, 52A}



\bigskip

\section{Introduction}

For every $K$ belonging to the class $\K^n$ of non-empty convex compact 
sets in $\R^n$, its quermassintegrals $W_i(K)$, for $i = 0,\dots, n$, 
are defined as the coefficients in the polynomial expansion
\begin{equation}{\mathcal H}^n (K + \rho B) = 
\sum_{i=0}^n \left({n \atop i} \right) W_i(K)\,\rho^i,
\end{equation}
where $\H^n$ denotes the Lebesgue measure on $\R^n$ and $K + \rho B$ 
is the Minkowski sum of $K$ plus $\rho$ times the unit Euclidean ball 
$B$. As special cases, $W_0$ is the Lebesgue measure $\H ^n$, $n W_1$ 
is the surface area, $2 \kappa_n^{-1} W _{n-1}$ is the mean width, 
and $\kappa_n^{-1} W_n = 1$ is the Euler characteristic 
(being $\kappa_n = \H^n(B)$).

\vskip2mm
The aim of this paper is to develope the notion of quermassintegrals 
for {\it quasi-concave} functions, as well as to enlighten their basic
properties. Quasi-concave functions $f$ on $\R ^n$ are defined by the inequality
$$
f((1 - \lambda) x_0 + \lambda x_1) \geq \min\{f(x_0),f(x_1)\},
\qquad \forall \, x_0,x_1 \in \R^n, \ \forall \lambda \in [0,1],
$$
and may also be described via the property that their level sets
$\{f \geq t\} = \{x \in \R^n: f(x)\geq t\}$
are convex. More precisely, we will work in the following class:
$$
\QC = \Big\{f :\R ^n \to [0, + \infty] :\ f\not\equiv 0\, , \ f
\hbox{ is quasi-concave, upper semi-continuous,}  \ 
\lim_{\|x\| \to + \infty} f (x) = 0\Big \},
$$
and also on the subclasses $\QC_\alpha$ of $\QC$ given by $\alpha$-concave functions, for $\alpha \in [- \infty, + \infty]$
(see Section 2.4 for details). 
The class $\QC$ can be considered a natural functional counterpart of $\K^n$:
in particular, for any $K \in \K^n$, its characteristic function  $\chi_K$ lies in $\QC$.

When passing from sets to (integrable) functions,  the role of the volume functional 
is played by the integral with respect to the Lebesgue measure:
\begin{equation}
I(f) = \int_{\sR^n} f(x) \,dx .
\end{equation}
This quite intuitive assertion, inspired by the equality 
$I(\chi_K) = \H^n(K)$, is commonly agreed and is also confirmed 
by several functional counterparts of geometric inequalities for convex bodies, in which the volume functional $\H^n(K)$ is replaced by 
the integral functional $I(f)$. As a significant example, one may
indicate the celebrated Pr\'ekopa-Leindler inequality 
\cite{B-L,Le,Pr1,Pr2,Pr3} 
(see also \cite{BaBo1,BaBo2,BuFr} for recent related papers), or the functional form of Blaschke-Santal\'o inequality \cite{AKM,Ball1}.

Less obvious is how to give a functional notion of the quermassintegrals
$W_i$ for $i>0$. The goodness of such a notion should be evaluated 
through the possibility of exporting to the functional framework 
the more relevant properties enjoyed by the quermassintegrals on $\K^n$.
The approach we propose goes exactly in this direction and relies on Cavalieri's principle: For every non-negative integrable function
$f$ on $\R ^n$,
$$
I(f) = \int_0^{+\infty} \H ^n \big(\{ f \geq t\}\big) \, dt.
$$
With a full consistency with the abstract Measure Theory (including 
its part dealing with integration over non-additive set functions), 
we define analogously the functionals
$$W_i(f) = \int_0^{+\infty} W_i\big(\{ f \geq t\}\big) \, dt, \qquad 
f \in \QC.
$$
The above definition is well-posed, since the mappings
$t \mapsto W_i \big (\{ f \geq t\} \big)$ are monotone increasing, 
as a consequence of the monotonicity of the functionals $W_i(\cdot)$ 
with respect to set inclusion. Actually, one can adopt the same natural
extension from sets to functions in more general situations: If $\Phi$ 
is any functional with values in $[0,+\infty)$, defined on $\K^n$ 
(or on the larger class of all Borel measurable subsets of $\R^n$), 
and if it is monotone increasing with respect to set inclusion, one 
can extend it to the class $\QC$ (respectively, to the class of 
all non-negative Borel measurable functions), by setting
\begin{equation}
\Phi (f) = \int_0^{+\infty} \Phi\big(\{f \geq t \}\big)\,dt.
\end{equation}

Definition (1.3) may look somewhat na\"ive if compared with 
previous notions existing in the literature for special quermassintegrals,
such as the perimeter or the mean width. These different definitions 
are rather based on the idea to mimic (1.1), by computing 
first order derivatives of the integral functional (1.2).
More precisely, starting from the the equalities
$$
\Per(K) = 
\lim_{\rho \to 0^+} \frac{\H ^n (K + \rho B) - \H ^n (K)}{\rho}\, , \qquad 
M(K) = \lim_{\rho \to 0^+} \frac{\H^n(B + \rho K) - \H^n(B)}{\rho}\, , 
$$
which are valid up to normalization constants for every $K \in \K^n$,
the following definitions have been considered in the recent works \cite{CoFr,Klartag-Milman05,Rotem1,Rotem2}, dealing especially with  log-concave functions:
$$
\Per(f) = \lim_{\rho \to 0^+} 
\frac{I(f \oplus \rho \cdot \varphi_n) - I(f)}{\rho}\ , \qquad 
M (f) = \lim_{\rho \to 0^+} 
\frac{I(\varphi_n \oplus \rho \cdot f) - I (f)}{\rho},
$$
where $\varphi_n$ denotes the density of the standard Gaussian measure 
on $\R^n$. Some more comments are in order to correctly understand 
the meaning of the above equalities. Firstly, the symbols $\cdot$ 
and $\oplus$ denote respectively a suitable multiplication by 
a nonnegative scalar and a suitable addition of functions, which can 
be defined so as to provide a natural extension of the 
usual Minkowski algebraic structure on $\K^n$ to functions, see Section 2 
for more details. 
Thus, the above definitions of perimeter and mean width, correspond to choose $\varphi_n$ as the functional counterpart of the unit ball on $\R^n$. Now, this choice may be somehow
disputable. To some extent, it is justified by the fact that the 
Gaussians turn out to be optimal in the functional version of 
meaningful geometric inequalities for which the Euclidean balls are 
optimal (see {\it e.g.} \cite{AKM}).

Notwithstanding, the investigation of the functional quermassintegrals introduced in (1.3) carried on in this paper, suggests a different 
point of view. As a starting point of this investigation, we consider, for a given $f \in \QC$ and any $\rho >0$, the functions
$$f _\rho (x) := \sup _{y \in B _\rho (x)}  f (y)\ ,
$$
where  $B _\rho (x)$ denotes the ball of radius $r$ centered at $x$. 
In fact, this is equivalent to perturb $f$ with the ``unit ball'' in the above mentioned algebraic structure, namely, if $f \in \QC_\alpha$,  it holds
$$
f _\rho = f \oplus \rho \cdot \Theta _\alpha(B) \,,
$$
being $\Theta _\alpha (B)$ the image of the unit ball through a natural isomorphic embedding of $\K ^n$ into $\QC _\alpha$. 
In particular, if $\alpha = - \infty$, meaning $f$ is merely quasi-concave, $\Theta _\alpha (B)$ is simply 
the characteristic function $\chi_B$. Therefore, in our perspective, $\chi _B$ is
the most natural functional equivalent of the ball $B$ in the class $\QC$. 
Actually, in Theorem 3.4, we prove that
a Steiner-type formula holds true 
for the mapping
\begin{equation}
\rho \mapsto I (f _\rho).
\end{equation}

More precisely,  we prove that such mapping is polynomial in $\rho$, 
and its coefficients are precisely the quermassintegrals defined in 
(1.3), see Theorem 3.4. In particular, up to normalization constants, the notions of 
perimeter and mean width of $f$ which are obtained from (1.3) with  $i=1$ and $i=n-1$,  correspond respectively 
to the coefficients of $\rho$ and of $\rho^{n-1}$ in the polynomial $I (f_\rho)$:
\begin{equation}
I(f _\rho) = 
I (f) + \Per(f)\,\rho + \dots + \frac{n \kappa _n}{2} M(f)\,\rho^{n-1} + \kappa_n(\max_{\sR^n} f)\, \rho^n\ .
\end{equation}
We then focus attention on the other main features of the 
quermassintegrals, dealing in particular with:

\medskip
-- concavity-like inequalities;

\medskip
-- integral-geometric formulae;

\medskip
-- valuation property;

\medskip
-- isoperimetric type inequalities.

\medskip
It is well-known that each of the functionals $W_i$'s satisfies on 
$\K^n$ the following Brunn-Minkowski type inequality:
\be
W_i((1 - \lambda) K_0 + \lambda K_1) \geq
\Big((1-\lambda) W_i(K_0)^{\frac{1}{n-i}} + \lambda W_i(K_1)^{\frac{1}{n-i}}\Big)^
{\frac{1}{n-i}}\quad \forall \,  K _0, K _1 \in \K ^n \,, \ \forall \lambda \in [0,1]\ .
\en

For short, this may be expressed as the property that the functional 
$\Phi = W_i$ is $\alpha$-concave on $\K^n$ with $\alpha = \frac{1}{n-i}$. 
For $i = 0$, namely for the Lebesgue measure, the functional 
counterpart of (1.7) is given by the dimension-free inequality due 
to Pr\'ekopa and Leindler and by its dimensional extension due to 
Brascamp and Lieb. We obtain a further generalization of these 
results (Theorems 4.2 and 4.7), which holds true for general 
monotone $\alpha$-concave functionals $\Phi$ extended from $\K^n$ 
to $\QC$ according to the formula (1.4). As a special case, 
we thus obtain Pr\'ekopa-Leindler-type inequalities for the 
functional quermassintegrals introduced in (1.3). 
On the example of the surface area, {\it i.e.} for the functional 
$\Phi = W_1$, the possibility of such generalization was already 
demonstrated in \cite{Bobkov}. As further examples of functionals satisfying a Brunn-Minkowski type inequality, let us mention the $p$-capacity 
of convex bodies in $\R ^n$ for 
$1 \leq p < n$ (with $\alpha = \frac{1}{n-p}$, see \cite{Bor3, CoSa}), the first non-trivial 
eigenvalue of the Laplacian with the Dirichlet boundary condition
(with $\alpha = -2$, see \cite{B-L}) and other similar functionals (see for instance \cite{Colesanti} and \cite{Salani}). 
These results link the study of quasiconcave functions to the theory of elliptic PDE's;  an example of the interaction between these 
subjects, particularly related to the matter treated here, can be found in \cite{Longinetti-Salani}.

Let us point out that our approach in order to prove Theorems 4.2 and 4.7 does not use induction on the 
dimension (nor mass transportation) as in the more typical proof of 
Pr\'ekopa-Leindler inequality, but is rather based on a new one-dimensional 
variant of it, inspired by a previous observation due to Ball \cite{Ball1}.
It is also remarkable that, as we show by constructing suitable counterexamples, this kind of  concavity property turns out to fail,
if one defines the perimeter of a function along the different line 
sketched above, namely as the derivative of the volume functional 
under Gaussian-type perturbations.

For what concerns integral-geometric results, we show that the 
Cauchy-Kubota formula for the quermassintegrals on $\K ^n$ can be 
suitably extended on $\QC$ (see Theorem 5.3). To that aim, we 
exploit as a crucial tool the concept of the functional projection 
introduced in \cite{Klartag-Milman05}. By combining it with definition 
(1.3), the desired extension turns out to be quite 
straightforward. To the best of our knowledge, this is the first step 
moved in bringing integral-geometric properties of convex bodies into 
a functional framework.

One of the most important characterizations of quermassintegrals is given 
by the celebrated Hadwiger's Theorem, which asserts that they generate 
the space of rigid motion invariant valuations on $\K^n$ which are 
continuous with respect to the Hausdorff metric (see \cite{Schneider}). The valuation property can be transferred in a natural 
way from sets to functions (replacing union and intersection by $\max$ 
and $\min$ operations, respectively, see Section 5 for 
details).  In Section 5 we check that the functionals defined in 
(1.3) are in fact valuations on $\mathcal Q^n$. Let us mention that recently some 
characterizations of valuations in various function spaces have been 
found,  see for instance \cite{Ludwig, Wright}. 

Besides concavity inequalities, and partly as a consequence of them, quermassintergrals verify various
inequalities of isoperimetric type; hence, having introduced a similar notion for functions, it is natural to ask
for corresponding results in the functional setting. In Section 6 we derive two possible versions of the standard
isoperimetric inequalities for quermassintegrals of quasi-concave and log-concave functions (see Theorems 6.1 and 6.2) along with a 
functional version of the Urysohn's inequality (Corollary 6.3).

\medskip
The outline of the paper is as follows.
After collecting some background material in Section 2, in Section 3 
we set and discuss our notion of functional quermassintegrals, and 
prove the corresponding Steiner formula. In Section 4 we deal with
generalized Pr\'ekopa-Leindler inequalities, while Section 5 is 
devoted to the integral-geometric formulae and the valuation property 
for functional quermassintergrals. Section 6 contains some concluding 
remarks on further properties related to isoperimetric and functional
inequalities.

\medskip

When this paper was in the final part of its preparation we learned by L. Rotem about the paper
\cite{Milman-Rotem}, where the authors present ideas and results, found independently, 
which partially overlap with those of the present paper. 

\medskip

\noindent{\bf Acknowledgment.} We wish to thank Paolo Salani for several discussions on the theme of quasi-concave functions, which 
gave a strong impulse to some of the ideas contained in this paper. 

\section{Preliminaries}

We work in the $n$-dimensional Euclidean space $\R^n$, $n\ge 1$, 
equipped with the usual Euclidean norm $\|\cdot\|$ and scalar product $(\cdot,\cdot)$. For $x\in\R^n$ and $r>0$, we set 
$B_r(x) = B(x,r) = \{y \in \R^n\, : \,\|y - x\| \leq r\}$, and
$B = B_1(0)$. We denote by $\relint(E)$ and $\cl(E)$ the relative 
interior and the closure of a set $E\subset\R^n$ respectively.

The unit sphere in $\R^n$ will be denoted by $\mathbb S^{n-1}$. 
For $k = 0,1,\dots,n$, $\H^k$ stands for the $k$-dimensional Hausdorff 
measure on $\R^n$. In particular, $\H^n$ denotes the usual Lebesgue 
measure on $\R^n$. 

\subsection{Convex bodies}
We denote by $\K ^n$ the class of all non-empty convex compact sets 
in $\R^n$ (called convex bodies). For the general theory of convex bodies, 
we refer the interested reader to the monograph \cite{Schneider}.

For every $K \in \K^n$, we denote by $\chi_K$ and $I_K$ respectively 
its characteristic and indicatrix functions, namely:
$$
\chi_K(x)=
\left\{
\begin{array}{lll}
\mbox{$1$, \ if $x\in K$,}\\
\mbox{$0$, \ if $x\notin K$,}
\end{array}
\right.
\quad
I_K(x)=
\left\{
\begin{array}{lll}
\mbox{$0$, \ \ \ if $x\in K$,}\\
\mbox{$+\infty$, \ if $x\notin K$.}
\end{array}
\right.
$$
Note that $I_K$ is convex. We will also use the notion 
of support function $h_K$ of a convex body $K$, defined by 
$$
h_K(x) = \sup_{y\in K}\, (x,y)\,.
$$
The class $\K^n$ is endowed with the algebraic structure based on 
the Minkowski addition. For $K$ and $L$ in $\K^n$, we set
$$
K+L = \{x+y\,|\,x\in K\,,\,y\in L\},
$$
while for $\lambda \ge 0$ and $K\in \K^n$, we set
$$
\lambda K = \{\lambda x\,|\,x\in K\}\,.
$$
It is worth noticing the following property connecting the Minkowski 
addition and support functions: For every $K,L \in \K^n$, 
and for every $\alpha,\beta \ge 0$,
$$
h_{\alpha K + \beta L} = \alpha h_K + \beta h_L\,.
$$

$\K^n$ can be endowed with the Hausdorff metric. The Hausdorff distance
between two convex bodies $K$ and $L$ can be simply defined as 
$$
\delta(K,L) = \|h_K-h_L\|_{L^\infty(\sfe)}
$$
(see \cite[Sec. 1.8]{Schneider}).

\subsection{Quermassintegrals of convex bodies}
In this subsection we collect basic properties and relations satisfied 
by the quermassintegrals. Recall that, for every $K \in \K ^n$, 
the quermassintegrals  $W_i(K)$, $i = 0, \dots, n$, represent 
the corresponding coefficients in the polynomial expansion (1.1). 
In particular, $W_0(K)=\H^n(K)$ is the volume of $K$, $W_n(K) = \kappa_n:= \H^n(B)$,
$n W_1(K) = \H^{n-1}(\partial K)$ is the surface area of $K$, and
$2\kappa_n^{-1} W_{n-1}(K)$ is the mean width, which is given by
$$
\int_{\sfe} (h_K(u)+h_K(-u))\,d\H^{n-1}(u).
$$

The quermassintegrals are invariant under rigid motions and continuous 
with respect to the Hausdorff distance. They also obey to the following
remarkable properties 
(where $K,\ K _0$ and $K _1$ denote arbitrary convex bodies in $\K ^n$).
\begin{itemize}
\item[(i)] {\it Homogeneity}.
$$
W_i(\lambda K) = \lambda^{n-i}\,  W_i(K)\quad\forall\ \lambda\ge0\,.
$$
\item[(ii)] {\it Monotonicity}.
$$
K _0 \subseteq K _1 \ \Rightarrow\ W_i(K_0) \leq W_i(K_1).
$$
\item[(iii)] 
{\it Brunn-Minkowski-type inequality}. For every $\lambda\in[0,1]$,
\be
W_i((1 - \lambda) K_0 + \lambda K_1) \geq
\Big((1-\lambda) W_i(K_0)^{1/(n-i)} + \lambda W_i(K_1)^{1/(n-i)}\Big)^{n-i}.
\en
Equivalently, the map
$\lambda \rightarrow W_i((1 - \lambda) K_0 + \lambda K_1)^\alpha$
is concave on $[0,1]$, where $\alpha = \frac{1}{n-i}$.
We will refer to this property as the $\alpha$-concavity of $W_i$. 
Note that in each case, $\alpha$ represents the reciprocal of the 
homogeneity order of the relevant quermassintegral.
The usual Brunn-Minkowski inequality corresponds to the case $i=0$. 

\item[(iv)] 
{\it Cauchy-Kubota integral formulae}. 
Given $k\in\{1,\dots,n-1\}$, let $\mathcal L^n_k$ be the set of all 
linear subspaces of $\R^n$ of dimension $k$, and let $dL_k$ denote the
integration with respect to the standard invariant probability measure 
on $\mathcal L^n_k$. Then, for every $i = 1,\dots,k$, we have
\begin{equation}
W_i(K)=c(i,k,n)\,\int_{\mathcal L^n_k} W_i(K|L_k)\,dL_k
\end{equation}
with a suitable constant $c(i,k,n)$. Here $K|L_k$ denotes the orthogonal projection of $K$ onto $L_k\in\mathcal \L^n_k$. An exhaustive 
presentation of these formulas (along with an explicit expression of 
the constant $c(i,k,n)$) may be found for instance in 
\cite{Schneider-Weil}. In the particular case $i=k=1$ we have 
the Cauchy integral formula for the perimeter:
$$
W_1(K) = c \int_{\mathbb S^{n-1}} \H^{n-1}(K|u^\perp)\, du\,,
$$
where $c$ is a constant depending on $n$ and $du$ indicates integration with respect to the invariant probability measure on the unit sphere.

\item[(v)] 
{\it Valuation property}. Every quermassintegral is a valuation on 
$\K^n$, {\it i.e.}, if $K_0$ and $K_1$ belong to $\K^n$ and are such that
$K_0\cup K_1\in\K^n$, then
\begin{equation}
W_i (K_0) + W_i (K _1) = W_i ( K _0 \cup K _1) + W_i ( K _0 \cap K _1) .
\end{equation}
According to a celebrated theorem by Hadwiger, this additivity property together with rigid motion invariance and continuity with respect to the Hausdorff distance (or monotonicity), characterizes linear combinations 
of quermassintegrals; see, for instance, 
Theorems 4.2.6 and 4.2.7 in \cite{Schneider}.
\end{itemize}

\subsection{$M$-means and $\alpha$-concave functions}
In order to introduce the class of $\alpha$-concave functions, we start 
with the definition of $\alpha$-means. Given 
$\alpha \in (-\infty,+\infty)$ and $s,t>0$, for every $u,v>0$ 
we first define
\begin{equation}
M_\alpha^{(s,t)}(u,v) :=
\left\{
\begin{array}{lll}
\mbox{$(s u^\alpha + t v^\alpha)^{1/\alpha}$, \ if $\alpha\ne0$,}\\
\mbox{\hskip7mm $u^sv^t$, \hskip12mm if $\alpha=0$.}
\end{array}
\right.
\end{equation}
For $\alpha\ge0$, definition (2.4) extends to the case when 
at least one of $u$ and $v$ is zero. If $\alpha<0$ and $uv=0$ 
(with $u,v\ge0$), we set $M_\alpha^{(s,t)}(u,v)=0$. 
In the extreme cases $\alpha = \pm \infty$, we set
$$
M_{-\infty}^{(s,t)}(u,v) := \min(u,v), \qquad
M_{+\infty}^{(s,t)}(u,v) := \max(u,v).
$$
The functions $u \rightarrow M_\alpha^{(s,t)}(u,v)$ and
$v \rightarrow M_\alpha^{(s,t)}(u,v)$ are non-decreasing.
If $u = + \infty$ or $v = +\infty$, the value $M_\alpha^{(s,t)}(u,v)$
is defined so that the monotonicity property is preserved. In particular,
$M_\alpha^{(s,t)}(+\infty,v) = M_\alpha^{(s,t)}(u,+\infty) =
+ \infty$ for every $v$ (including $v=+\infty$) in case $\alpha > 0$.
We also put 
$M_\alpha^{(s,t)}(+\infty,0) := M_\alpha^{(s,t)}(0,+\infty) = 0$ 
for $\alpha\le0$.

\vskip2mm
The $\alpha$-mean of $u,v\ge0$, with weight $\lambda\in(0,1)$ is 
defined as
$$
M_\alpha^{(\lambda)}(u,v) = M_\alpha^{(1-\lambda,\lambda)}(u,v)\,.
$$
The particular cases $\alpha=1,0,-1$ correspond to the arithmetic, 
geometric and harmonic mean, respectively. In general,
the functions $\alpha \rightarrow M_\alpha^{(\lambda)}(u,v)$
are non-decreasing. Note, however, that this property fails for the
functions $\alpha \rightarrow M_\alpha^{(s,t)}(u,v)$ with $s+t \neq 1$.

\bigskip

For $\alpha\in[-\infty,+\infty]$, we denote by $\mathcal C_\alpha$ the 
family of all functions $f: \R^n \to [0, + \infty]$ 
which are not identically zero and are $\alpha$-concave, meaning that  
$$
f((1 - \lambda) x + \lambda y) \geq M_\alpha^{(\lambda)}(f(x),f(y)), 
\qquad \forall \, x,y \ \, {\rm such \ that} \  f(x)f(y)>0,\ \ \forall \,
\lambda \in (0,1).
$$

\bigskip
The same definition may be given when $f$ is defined on a convex 
subset of $\R^n$. Note that, as a straightforward consequence of the
monotonicity property of the $\alpha$-means with respect to $\alpha$, 
we have $\mathcal C_\alpha \subseteq \mathcal C_{\alpha'}$ if 
 $\alpha'\le\alpha$. 

The following particular cases of $\alpha$ describe canonical classes 
of $\alpha$-concave functions:
\begin{itemize}
\item[] 
$\mathcal C_{- \infty}$  is the largest class of quasi-concave functions;
\item[] 
$\mathcal C_{0}$  is the class of log-concave functions;
\item[] 
$\mathcal C_{1}$  is the class of concave functions on convex sets $\Omega$ 
(extended by zero outside $\Omega$);
\item[] 
$\mathcal C_{+ \infty}$  is the class of multiples of characteristic
functions of convex sets $\Omega \subset \R^n$. 
\end{itemize}

Any function $f\in \mathcal C_\alpha$ is supported on the (nonempty) convex set
$K_f = \{f>0\}$, and if $\alpha > -\infty$, it is continuous in the 
relative interior $\Omega_f$ of $K_f$. If $\alpha$ is finite and nonzero,
it has the form $f = V^{1/\alpha}$, where $V$ is concave on $\Omega _f$ in case 
$\alpha > 0$, and is convex in case $\alpha < 0$; for $\alpha = 0$, 
the general form is $f = e^{-V}$ for some convex function $V$ on $\Omega_f$.


\subsection{Algebraic structure of the class of $\alpha$-convex functions}

For any $\alpha \in [ - \infty, + \infty]$, we are going to introduce 
in $\mathcal C_\alpha$
an addition and a multiplication 
by positive reals, which extend the usual Minkowski algebraic structure on $\K^n$. 

Let be given $f,g \in\mathcal C_\alpha$ and $s,t > 0$.  
If  $\alpha\le 0$, we
put
\be
(s\cdot f\oplus t\cdot g)(z): = \sup\left\{M_\alpha^{(s,t)}(f(x),g(y)): \,
z = s x + t y\right\};
\en
if $\alpha>0$, we put
\be
(s\cdot f\oplus t\cdot g)(z): = 
\left\{
\begin{array}{lll}
\sup\left\{M_\alpha^{(s,t)}(f(x),g(y)): \,
z = s x + t y, \ f(x)g(y)>0\right\} \mbox{ \ \ if \, $z\in sK_f+tK_g$,}
\\
\mbox{$0$ \hskip90mm otherwise.}
\end{array}
\right.
\en
Note that (2.6) is also applicable in case
$\alpha\le0$, since $M^{(s,t)}_\alpha(u,v)=0$, whenever $uv=0$;
in this sense (2.6) is more general than (2.5). 

Clearly the operations $\oplus$ and $\cdot$ depend on $\alpha$. However for simplicity we will not indicate this 
dependence explicitly, unless it is strictly needed. In particular, this abuse of notation is consistent with the following immediate relation: For all non-empty sets $K$ and $L$ in $\R^n$ and all $s,t > 0$,
$$
s\cdot \chi_K\oplus t\cdot \chi_L = \chi_{sK + tL}
$$
(in particular, in this case the left-hand side does not depend on $\alpha$).


The operations $\oplus$ and $\cdot$ may also be used for arbitrary 
non-negative, not identically zero functions, without any convexity 
assumption. For any fixed $\alpha \in [ - \infty, + \infty]$, they are easily checked to
enjoy the following general properties: 

\vskip2mm
\begin{itemize}
\item[(i)] {\it Commutativity}. 
$s\cdot f\oplus t\cdot g = t\cdot g\oplus s\cdot f$.
\item[(ii)] {\it Associativity}. \ 
$(s \cdot f \oplus t\cdot g) \oplus u\cdot h = 
s \cdot f \oplus (t\cdot g \oplus u\cdot h)$.
\item[(iii)] {\it Homogeneity}. \ 
$s \cdot f \oplus t\cdot g = (s+t)^{1/\alpha}\,
\big(\frac{s}{s+t} 
\cdot f \oplus \frac{t}{s+t} \cdot g\big)$ \ \ $(\alpha \neq 0)$.
\item[(iv)] {\it Measurability}. \ 
$s \cdot f \oplus t\cdot g$ is Lebesgue measurable as long as $f$ and $g$ are Borel measurable. 

\end{itemize}

Next, we show that every class $\mathcal C_\alpha$
is closed under the introduced operations.

\vskip5mm
\begin{prop}
If $f,g\in\mathcal C_\alpha$ and $s,t > 0$, 
then $s\cdot f\oplus t\cdot g\in\mathcal C_\alpha$.
\end{prop}

\proof First let $\alpha$ be non-zero. Using the homogeneity property 
(iii), it suffices to consider the case $s+t=1$. We set for brevity 
$$
u(x,y) = M_\alpha^{(s,t)}(f(x),g(y)), \qquad x \in K_f, \ y \in K_g,
$$
and, for $z \in K = sK_f+tK_g$, let
$$
h(z) = (s\cdot f\oplus t\cdot g)(z) = \sup\left\{u(x,y): \,
z = sx + ty, \ x \in K_f, y \in K_g\right\},
$$
putting $h = 0$ outside $K$.

We claim that the function $u$ is $\alpha$-concave on the convex
supporting set $K_f \times K_g$. Indeed, if additionally $\alpha$ 
is finite, taking
$(x,y) = s'(x_1,y_1) + t'(x_2,y_2)$ with $s',t' > 0$, $s'+t' = 1$
and $x_1, x_2 \in K_f$, $y_1, y_2 \in K_g$, we have
\bee
u(x,y) & = &
M_\alpha^{(s,t)}\left(f(s' x_1 + t' x_2),g(s' y_1 + t' y_2)\right) \\
 & \geq &
M_\alpha^{(s,t)}\left(M_\alpha^{(s',t')}(f(x_1),f(x_2)),
M_\alpha^{(s',t')}(g(y_1),g(y_2))\right) \\
 & = &
\Big(s
\left(s' f(x_1)^\alpha + t' f(x_2)^\alpha\right) +
t
\left(s' g(y_1)^\alpha + t' g(y_2)^\alpha\right)\Big)^{1/\alpha} \\
 & = &
\Big(s'\left(sf(x_1)^\alpha + t g(y_1)^\alpha\right) +
t'\left(sf(x_2)^\alpha + t g(y_2)^\alpha\right)\Big)^{1/\alpha}
\\
 & = &
M_\alpha^{(s',t')}\left(M_\alpha^{(s,t)}(f(x_1),g(y_1)),
M_\alpha^{(s,t)}(f(x_2),g(y_2))\right) \\
 & = &
M_\alpha^{(s',t')}(u(x_1,y_1),u(x_2,y_2)).
\ene
Thus,
$$
u(s'(x_1,y_1) + t'(x_2,y_2)) \geq 
M_\alpha^{(s',t')}(u(x_1,y_1),u(x_2,y_2)),
$$
which means $\alpha$-concavity of $u$ on $\R^{2n}$ 
(if we define it to be zero outside $K_f \times K_g$).

With corresponding modifications, or using continuity and monotonicity
of the function $M_\alpha$ with respect to $\alpha$,
we have a similar property of the function $u$ in the remaining cases.

Now, for $z \in K$, fix a decomposition $z = sz_1 + t z_2$, 
$z_1,z_2 \in K$. Using truncation, if necessary, we may assume that both 
$f$ and $g$ are bounded, so that $h$ is bounded, as well. Then, 
given $\ep>0$, choose $x_1,x_2 \in K_f$, $y_1,y_2 \in K_g$ such that
$z_1 =s x_1 + t y_1$,
$z_2 =s x_2 + t y_2$, and
$$
h(z_1) \leq u(x_1,y_1) + \ep, \qquad h(z_2) \leq u(x_2,y_2) + \ep.
$$
Since the function $u$ is $\alpha$-concave, setting $x = s x_1 + t x_2$ and
$y = s y_1 + t y_2$, we get
$$
u(x,y)
  \geq
M_\alpha^{(s,t)}(u(x_1,y_1),u(x_2,y_2))
  \geq
M_\alpha^{(s,t)}\left((h(z_1)- \ep)^+,(h(z_2)- \ep)^+\right)\, .
$$
Letting $\ep \rightarrow 0$, the latter yields
$$
u(x,y) \geq M_\alpha^{(s,t)}(h(z_1),h(z_2))\,.
$$
It remains to note that $sx + ty = s z_1 + t z_2 = z$,
which implies $u(x,y) \leq h(z)$. 

Now, let $\alpha = 0$, in which case we should work with
$$
u(x,y) = f(x)^s g(y)^t, \qquad x,y \in \R^n,
$$
and with a similarly defined function $h$. Again, for 
$(x,y) = s'(x_1,y_1) + t'(x_2,y_2)$, we have, using the log-concavity
of $f$ and $g$,
\bee
u(x,y) & = &
f(s' x_1 + t' x_2)^s \, g(s' y_1 + t' y_2)^t \\
 & \geq &
f(x_1)^{s s'} f(x_2)^{st'} \, g(y_1)^{ts'} g(y_2)^{tt'}
 \ = \
M_0^{(s',t')}(u(x_1,y_1),u(x_2,y_2)).
\ene
This means that $u$ is log-concave on $\R^{2n}$. The rest
of the proof is similar to the basic case. \qed

\bigskip

In the next remarks we collect further comments on the operations $\oplus$ and $\cdot$, more specifically on their relationship with the usual Minkowski structure in $\K ^n$, and on their interpretation in the two
special cases $\alpha = - \infty$ and $\alpha = 0$. 

\begin{remark} {\rm Equipped with quermassintegral in (2.5)-(2.6), and in view of Proposition 2.1, 
$\mathcal C_\alpha$ can be seen as an extension of 
$\K^n$ which preserves its algebraic structure. More precisely,
the mappings $\Theta _\alpha : \K ^n \to {\mathcal C} _\alpha$ defined by
$$\Theta _\alpha (K) := \begin{cases} e ^ {- I _K} & \hbox{ if } \alpha = 0 \\ \noalign{\smallskip}
I _K ^ {-1} & \hbox{ if } \alpha \neq 0 
\end{cases}
$$
are isomorphic embeddings of $\K ^n$ (endowed with the Minkowski structure) into ${\mathcal C} _\alpha$ 
(endowed with the operations $\oplus$ and $\cdot$). }
\end{remark}


\begin{remark} {\rm In $\mathcal C_{-\infty}$, quermassintegral in (2.5) can be characterized through the Minkowski 
addition of the level sets $K_f(r) = \{x\in\R^n\,:\,f(x)> r\}$.
Namely, for $f,g\in\mathcal C_{-\infty}$ and $s,t>0$, the functional 
equality
$$
h(z) = (s\cdot f\oplus t\cdot g)(z) = \sup\{\min\{f(x),g(y)\}:sx+ty=z\}
$$
is equivalent to the family of set equalities
\begin{equation}
K_h(r) = sK_f(r) + tK_g(r) \qquad \forall \, r >0\,.
\end{equation}
Note that for a general value of $\alpha$, 
we only have the following set inclusion, valid if $s+ t =1$:
\begin{equation}
K_h(r) \supset sK_f(r) + tK_g(r) \qquad  \forall \, r >0\,.
\end{equation}
}
\end{remark}

\bigskip
\begin{remark} {\rm In $\mathcal C_0$,  the operation $\oplus$ 
(defined as in (2.5) with $t=s=1$) is related to the operation introduced in 1991 
by Maurey. More precisely, starting with 
$U,V:\R^n \rightarrow (-\infty,+\infty]$, we get
\begin{equation}
e^{-U} \oplus e^{-V} = e^{-W},
\end{equation}
where
$$
W(z) = \inf_x \ [U(z-x) + V(x)]
$$
represents the infimum-convolution of $U$ and $V$. If these functions
are convex, so is $W$ (as we also know from Proposition 2.1). 
This fact is crucial in the study of the so-called ``convex"
concentration for product measures, {\it cf.} \cite{Maurey}.
}
\end{remark}

\subsection{Pr\'ekopa-Leindler and Brascamp-Lieb Theorems}
The following well-known result due to Pr\'ekopa and Leindler
\cite{Le,Pr1,Pr2,Pr3} is a functional extension of the classical Brunn-Minkowski inequality.

\vskip5mm
\begin{teo} Let $\lambda \in (0,1)$. Let  $f,g,h$ be 
non-negative measurable functions on $\R^n$. If 
\be
h(M_1^{(\lambda)}(x,y)) \geq  M_0^{(\lambda)}(f(x),g(y)) \qquad \forall x, y \in \R ^n\, ,
\en
then
\be
\int h \geq M_0^{(\lambda)}\bigg(\int f, \int g \bigg).
\en
\end{teo}


Given non-empty Borel sets 
$A,B \subset \R^n$, and $\lambda \in (0,1)$, 
by applying the above result with $f = \chi_A$, $g = \chi_B$, and $h = \chi_{(1-\lambda) A + \lambda B}$
(after noticing that $h$ is Lebesgue measurable), 
one gets
\be
\H^n ((1-\lambda) A + \lambda B) \geq \H^n(A)^{1-\lambda} \H^n(B)^\lambda.
\en
This is a multiplicative variant of the Brunn-Minkowski inequality
\be
\H^n((1-\lambda) A + \lambda B) \geq
\left((1-\lambda) \H^n(A)^{1/n} + \lambda\H^n(B)^{1/n}\right)^n
\en
with convexity parameter $\alpha = 1/n$ (which is optimal). Though 
in principle (2.12) is weaker (2.13), using the homogeneity of the 
volume it is easy to derive (2.13) from (2.12). However, the difference
between (2.13) and (2.12) suggests a different, dimension-dependent 
variant of Theorem 2.5, which would directly yield (2.13) when applied 
to characteristic functions. Such a variant is known and is recalled 
in Theorem 2.6 below. It was proposed by Brascamp and Lieb \cite{B-L} 
and somewhat implicitly in Borell \cite{Bor1,Bor2}; 
{\it cf.} also \cite{DG} and \cite{D-U}.

\begin{teo} 
Let $\lambda \in (0,1)$ and let $\alpha \in [-\frac{1}{n},+\infty]$. 
Let $f,g,h$ be non-negative measurable functions on $\R^n$. If
\be
h(M_1^{(\lambda)}(x,y)) \geq M_\alpha^{(\lambda)} (f(x),g(y)),
\qquad \forall \, x,y \,  \hbox{ such that} \  f (x)g(y) > 0,
\en
then
\be
\int h \geq M_\beta^{(\lambda)}\bigg(\int f, \int g \bigg) 
\qquad \hbox{where} \ \   \beta := \frac{\alpha}{1 + \alpha n}.
\en
In the extreme cases
$\alpha = -\frac{1}{n}$ and $\alpha = +\infty$, the definition of 
$\beta$ in $(2.15)$ is understood respectively as $\beta = -\infty$ 
and $\beta = \frac{1}{n}$.
\end{teo}

Since $\beta = 0$ for $\alpha = 0$, Theorem 2.6 includes Theorem 2.5 
as a particular case. Note also that, if $A, B$ and $\lambda$ are as above, by applying Theorem 2.6 with 
$\alpha = + \infty$, $f = \chi_A$, $g = \chi_B$ and 
$h = \chi_{(1-\lambda) A + \lambda B}$, one obtains directly 
the Brunn-Minkowski inequality in its dimension-dependent form (2.13).


We point out that, under additional assumptions on $f$ and $g$, the value of $\beta$
in (2.15) may be improved. For instance, in dimension $n=1$, if
${\rm ess\,sup}\, f(x) = {\rm ess\,sup}\, g(x) = 1$, then one may
take $\beta = 1$ regardless of $\alpha$, see for instance \cite{Bobkov-Ledoux}.
Without additional constraints, the value of $\beta$ in (2.15) is optimal.
For instance, for $n=1$ and $\alpha = 0$, take $f(x) = ae^{-x} \chi_{(0,+\infty)}(x)$ 
and $g(x) = be^{-x} \chi_{(0,+\infty)}(x)$, where $a$ and $b$ are positive parameters.
In this case, the function $h(x) := M_0^{(\lambda)}(a,b)\,e^{-x} \chi_{(0,+\infty)}(x)$
satisfies (2.10), and (2.11) becomes equality.

\bigskip
As a further natural generalization of Theorem 2.6, one can consider the case when $\lambda$ and $(1- \lambda)$
are replaced by arbitrary positive parameters $s$ and $t$, not necessarily satisfying the condition $s+t = 1$. 
Assume $\alpha \neq 0$, and $\alpha < + \infty$. If non-negative measurable functions  $f,g,h$ satisfy the inequality $
h(M_1^{(s,t)}(x,y)) \geq M_\alpha^{(s,t)} (f(x),g(y))$ for all $x, y$ such that $f (x) g (y) >0$, then the function 
$$
\tilde h(z) := \frac{1}{(s+t)^{1/\alpha}}\, h((s+t)\,z)
$$ 
is easily checked to satisfy the hypothesis (2.14) with 
$\lambda = \frac{t}{s+t}$. Hence, by applying Theorem 2.6, we arrive at
the following statement  (where also the case 
$\alpha = +\infty$ can be easily included as a limit):

\begin{teo} Let $s, t>0$ and let $\alpha \in [-\frac{1}{n},+\infty]$, 
$\alpha \neq 0$. Let $f,g,h$ be non-negative measurable functions on 
$\R^n$.  If
$$
h(M_1^{(s,t)}(x,y)) \geq M_\alpha^{(s,t)} (f(x),g(y)),
\qquad \forall \, x,y \,  \hbox{ such that} \  f (x)g(y) > 0,
$$
then
$$
\int h \geq M_\beta^{(s,t)}\bigg(\int f,\int g \bigg) 
\qquad \hbox{where} \ \   \beta := \frac{\alpha}{1 + \alpha n}.
$$
In the extreme cases
$\alpha = -\frac{1}{n}$ and $\alpha = +\infty$, the value of $\beta$ has to be understood as in Theorem 2.6.
\end{teo}

\bigskip

We observe that, using the operations $\oplus$ and $\cdot$ introduced in the previous section, 
Theorem 2.7 (and similarly also Theorems 2.5 and 2.6) can be written in a more compact form as the inequality
$$
\int \big(s \cdot f \oplus t \cdot g\big) \geq  
M_\beta^{(s,t)}\bigg(\int f, \int g \bigg), \qquad \hbox{ where } 
\alpha \in  [ - \frac{1}{n}, + \infty]\, ,\ \alpha \neq 0, \hbox{ and } 
\beta = \frac{\alpha}{1 + \alpha n}, 
$$
holding true for all non-negative Borel measurable functions $f$ and $g$ 
on $\R^n$, and for all $t,s > 0$ (the assumption $\alpha \neq 0$ 
may be removed when $t+s=1$).

In particular, taking $s=t=1$, and 
replacing first $f$, $g$ 
respectively with  $f^{1/\alpha}$, $g^{1/\alpha}$, and then  
$\alpha$ with $\frac{1}{\alpha}$, one gets the following inequality
\be
\int 
\Big ( \sup \{ f (x) + g (y) \ :\ x + y = z\, , f (x) g (y) >0 \} \Big ) ^ \alpha
 \geq  
\bigg[\Big(\int f^\alpha\Big)^{\frac{1}{\alpha + n}} +
\Big(\int g^\alpha\Big)^{\frac{1}{\alpha + n}}\bigg]^{\alpha + n},
\en
where  $\alpha \geq 0$  or $\alpha \leq -n$.

In dimension $n=1$ and for the range $\alpha > 0$, this inequality
was obtained in 1953 by Henstock and Macbeath as part of their proof
of the Brunn-Minkowski inequality, {\it cf.} \cite{Henstock-Macbeath}.
Indeed,  stated in $\R^n$ for characteristic 
functions $f = \chi_A$, $g = \chi_B$, and with $\alpha = 0$, 
(2.16) gives back 
$$
\H^n(A + B) \geq \left(\H^n(A)^{1/n} + \H^n(B)^{1/n}\right)^n.
$$


\section{Functional notion of quermassintegrals and Steiner-type formula}

Let us introduce the following class of admissible functions
$$
\QC = \Big\{f : \R ^n \to [0, + \infty] \ :\ f\not\equiv 0\,, \ f \ \hbox{is quasi-concave,  
upper semicontinuous}, \ \lim_{\|x\| \to +\infty} f(x)=0\Big\}.
$$

We also consider the subclasses formed by the functions in $\QC$ which are $\alpha$-concave: 
$$
\mathcal Q^n_\alpha =\QC \cap {\mathcal C} _\alpha\, , \qquad \alpha\in [ - \infty,  +\infty].
$$
In particular, $\QC=\QC_{-\infty}$.

Note that, if $f$ is quasi-concave, the property $\lim_{\|x\| \to +\infty} f(x)=0$ is necessary to keep 
$I(f)$ finite (we recall that $I(f)$ is just the integral of $f$ on $\R^n$).
Indeed, the vanishing of $f$ at infinity may be equivalently formulated as the boundedness of all the level sets $\{f\geq t\}$: if $I (f)$ is finite, then all such convex sets have finite Lebesgue measure and are therefore bounded.

We also observe that, if $f \in \QC$, the level sets $\{f \geq t\}$ are convex closed sets, because
$f$ is quasi-concave and upper semicontinuous; 
since $f$ is vanishing at infinity, these sets are also compact. Hence,
$\sup_x f(x)$ is attained at some point, and one may freely speak about 
the maximum value of $f$ (which in general may be finite or not).
In addition, all quermassintegrals of the sets $\{ f \geq t \}$ are well-defined
and finite, so that we are allowed to give the the following definition. 

\begin{definition}
{\rm Let $f \in \QC$. For every $i=0,\dots,n$, we define the 
{\it i-th quermassintegral of $f$} as
\begin{equation}
W_i(f) := \int_0^{+\infty} W_i \big(\{f \geq t\} \big)\,dt = 
\int_0^{+\infty} W_i\big(\cl\{f > t\}\big)\,dt.
\end{equation}
In particular,
$$ 
I (f) = W_0 (f) = \int_0^{+\infty} \H^n\big(\{f\geq t\}\big)\,dt.
$$
As further special cases, 
by analogy with convex bodies, we define the {\it perimeter}, the {\it mean width} and 
the {\it Euler characteristic} of $f \in \QC$ respectively as
$$
\begin{array}{ll} & 
\displaystyle{\Per(f) = n W_1 (f) = 
\int_0^{+\infty} \Per\big(\{ f\geq t\}\big)\,dt}, \cr
\noalign{\smallskip} & 
\displaystyle{M(f) = 2\kappa_n^{-1}\, W_{n-1}(f) = 
\int_0^{+\infty} M\big(\{f\geq t\}\big)\,dt}, \cr
\noalign{\smallskip} & 
\displaystyle{\chi(f) = \kappa_n^{-1}\, W _n(f) = \max_{x \in \sR^n} f(x)}.
\end{array}
$$}
\end{definition}

Let us emphasize that the two integrals in (3.1) do coincide, so that  we 
may use any of them at our convenience. To see this fact, one 
may use the inclusion ${\rm cl}\{f>t\} \subseteq \{f \geq t\}$, which 
ensures that the second integral in (3.1) is dominated by 
the first one (applying the monotonicity property of $W_i$). 
On the other hand, for any $\ep > 0$, we have
$\{f \geq t+\ep\} \subseteq \{f>t\} \subseteq {\rm cl}\{f>t\}$, which
yields
$$
\int_\ep^{+\infty} W_i \big(\{f \geq t\} \big)\,dt \leq 
\int_0^{+\infty} W_i\big(\cl\{f > t\}\big)\,dt.
$$
Letting $\ep \rightarrow 0$, we obtain that the first integral 
in (3.1) is dominated by the second one, as well.

\subsection{Basic properties}

Let us mention a few general properties of the functional
quermassintegrals, which follow immediately from Definition 3.1. 

\vskip2mm
\begin{itemize}
\item[(i)] {\it Positivity}. \
$0 \leq W_i(f) \leq +\infty$.
\item[(ii)] {\it Homogeneity under dilations}. \ 
$W_i(f_\lambda) = \lambda^{n-i}\, W_i(f)$, where
$f_\lambda(x) = f(x/\lambda)$, $\lambda>0$.
\item[(iii)] {\it Monotonicity}. \ 
$W_i(f) \leq W_i(g)$, whenever $f \leq g$.
\end{itemize}

For what concerns the finiteness of the quermassintergals, 
the problem of characterizing those functions in $\QC$ whose all quermassintegrals are finite
seems to be an interesting question.  Let us examine what happens in this respect within the subfamily of radial functions.

\begin{example}{\rm  Let $f \in \QC$ 
be a spherically invariant function. Equivalently, it has the form
$$
f(x) = F(|x|), \quad x \in \R^n,
$$
where $F:[0,+\infty) \rightarrow [0,\Lambda]$ is a non-increasing upper
semi-continuous function vanishing at infinity, with maximum 
$\Lambda = F(0)$, finite or not.

Incidentally, this example shows that quasi-concave functions do not 
need to be continuous on their domain, nor to be in $L ^1(\R^n)$, so that
it may be $I (f) = + \infty$. 

Define the inverse function $F^{-1}:(0,\Lambda] \rightarrow [0,+\infty)$
canonically by
$$
F^{-1}(t) = \min\{r>0: F(r) \geq t\}, \qquad 0 < t < \Lambda.
$$
Since $\{f \geq t\} = F^{-1}(t) B$, we have
$W_i(\{f \geq t\}) = \kappa_n \big(F^{-1}(t)\big)^{n-i}$.
Integrating this equality over $t$, we arrive at the formula
$$
W_i(f) = \kappa_n \int_0^{+\infty} r^{n-i}\,dF(r), \qquad
i = 0,1,\dots,n,
$$
where $F$ may be treated as an arbitrary positive measure on 
$(0,+\infty)$, finite on compact subsets of the positive half-axis. 
Hence, the quermassintegrals of the
function $f$ are described as the first $n$ moments of $F$
(up to the normalization constant $\kappa_n$).

In particular, we see that 
the finiteness of $W_n(f)$ is equivalent to the finiteness of the measure $F$ (namely to the condition $\Lambda < +\infty$),
whereas the finiteness of $W_0(f)$
is equivalent to $\int_0^{+\infty} r^n\,dF(r) < +\infty$.
Thus we can conclude that the quermassintegrals $W_i(f)$ are finite for all $i = 0, \dots,  n$, 
if and only if they are finite for $i=0$ and $i=n$.
} 
\end{example}

\vskip2mm
The above example suggests a simple way to find upper bounds on the quermassintegrals in the general case. Namely, 
the monotonicity property (iii) stated above readily yields:

\begin{prop} {\it Given a function $f \in \QC$, define
$\mu _f(r) = \max_{\|x\| \geq r} f(x)$, $r > 0$. Then
$$
W_i(f) \leq \kappa_n \int_0^{+\infty} r^{n-i}\,d\mu _f(r), \qquad
i = 0,1,\dots,n.
$$
In particular, all quermassintegrals of $f$ are finite, provided $f$ 
is bounded and $\int_0^{+\infty} r^n\,d\mu _f(r) < +\infty$.
}
\end{prop}

\bigskip
\subsection{Steiner formula}

Let $f \in \QC$. For $\rho >0$, consider the function  
$$
f_\rho (x) = \sup_{y \in B_\rho(x)} f(y)\, .
$$
If $f \in \QC _\alpha$, using the operations $\oplus$ and $\cdot$ introduced in Section 2.4 on the class $\mathcal C_\alpha$, and the isomorphic embeddings $\Theta _\alpha$ of Remark 2.2, the  function $f _\rho$ may also be rewritten as
$$
f _\rho = f \oplus \rho \cdot \Theta _\alpha(B) \,
$$
(recall that $B = B _1 (0)$, $\Theta _0 (B) = \chi _B$, and $\Theta _\alpha (B) = I _B ^ {-1}$ for $\alpha \neq 0)$. 
Therefore, the function $f _\rho$ can be seen as a perturbation of $f$ through the unit ball. 
Actually, the next result provides a functional analogue of the Steiner formula, stating that 
the integral of $f_\rho$ admits 
a polynomial expansion in $\rho$, with coefficients given precisely 
by the functional quermassintegrals $W_i(f)$'s.

\begin{teo}{\rm (Steiner-type formula)}
Let $f \in \QC$. For every $\rho>0$, there holds
\begin{equation}
I(f_\rho) = \sum_{i=0}^n \left({n \atop i}\right) W_i(f)\,\rho^i.
\end{equation} 
\end{teo}

Before giving the proof of Theorem 3.4, let us point out that, as a consequence of  (3.2), the following properties turn out to be equivalent to each other:
\begin{itemize}
\item[(i)]  $W_i (f)< + \infty$ \ $\forall i = 0, \dots, n $;
\item[(ii)] $I ( f _\rho) < + \infty$ for some $\rho>0$;
\item[(iii)] $I ( f _\rho) < + \infty$ for all $\rho>0$. 
 \end{itemize}

In particular, the condition $I (f) < + \infty$ is not sufficient to guarantee that  
$I(f_\rho) < +\infty$ (as the latter condition implies the boundedness of $f$). 
A simple sufficient condition is  for instance that 
$f$
is of class $C^1(\R^n)$, with $I (f) < + \infty$ and
$$
\int_{\R^n} \max_{y \in B_\rho(x)} \|\nabla f(y)\|\,dx < +\infty\, ;
$$
indeed, by using the inequality 
$f_\rho(x) \leq 
f(x) + \max_{y \in B_\rho(x)} \|\nabla f(y)\|$, it follows that 
$I(f_\rho) < +\infty$.

Whenever $I(f_\rho)$ is finite, as an immediate consequence of 
Theorem 3.4, the quermassintegrals $W_i (f)$ can be expressed through differential  formulae involving $I(f_\rho)$. In particular, it holds
\begin{equation}
\Per(f) = \lim_{\rho \to 0^+} \frac{I(f_\rho) - I(f)}{\rho}
\end{equation}
and
\begin{equation}
M(f) = \frac{2}{n \kappa_n} \lim_{\rho \to +\infty} 
\frac{I(f_\rho) - (\kappa_n \max_{\sR^n} f)\, \rho^n}{\rho^{n-1}}.
\end{equation}

\bigskip
\begin{remark}
{\rm Let $f \in \QC$. Denote by $K_f$ the support set $\{f>0\}$, 
by $|Df|(\R ^n)$ the total variation of $f$ as a $BV$ function on $\R^n$, 
and by $f_+$ the interior trace of $f$ on $\partial K_f$. Then
\begin{equation}
\Per(f) = \int_0^{+\infty} \Per(\{f \geq t \})\,dt = 
|Df|(\R^n) = \int_{K_f} |\nabla f|\,dx + \int_{\partial K_f} f_+ \,d\H^{n-1},
\end{equation}
where we have used the definition of $\Per(f)$ and the coarea formula. 
This formula is simplified to
$$
\Per(f) = \int_{\R^n} |\nabla f|\,dx,
$$
if $f$ is continuously differentiable on the whole $\R^n$ (which also 
follows from (3.3) in case $I(f_\rho) < +\infty$, for some $\rho>0$).
We point out that (3.5) may be seen as a variant of the integral representation formula given by Theorem 4.6 in \cite{CoFr}: in fact, 
(3.5) can be derived ``formally''  by applying Theorem 4.6 in \cite{CoFr} beyond its assumptions (more precisely, by taking therein $\psi (y) = |y|$).}
\end{remark}

\bigskip
{\it Proof of Theorem 3.4}. We start from the well-known elementary
identity (which is often used in derivation of various Sobolev-type 
inequalities)
\begin{equation}
\big\{f_\rho > t\big\} = \big\{f > t\big \} + \rho B \qquad (\rho,t>0).
\end{equation}

Define the sets
$$
\Omega^t = \{f>t\}, \qquad \Omega_\rho^t = \{f_\rho > t\}, \qquad 
K^t = \cl \Omega^t, \qquad K_\rho^t = \cl \Omega_\rho^t.
$$ 
Since $f \in \QC$, the convex sets $\Omega^t$ are bounded, so are 
$\Omega_\rho^t$, and one has $\H^n(\Omega_\rho ^t) = \H^n(K_\rho ^t)$. 
Then, by virtue of Cavalieri's principle, we can express $I(f_\rho)$ as
$$
I(f_\rho) = \int_0^{+\infty} {\mathcal H}^n (\Omega_\rho^t) \, dt = 
\int_0^{+\infty} {\mathcal H}^n (K_\rho^t) \, dt.
$$ 
By (3.6), we have
$$
K_\rho ^t = \cl (\Omega _\rho^t) = \cl (\Omega^t + \rho B) =  
\cl (\Omega^t) + \rho B =  K^t + \rho B. 
$$
Hence,
$$
I(f_\rho) = \int_0^{+\infty} {\mathcal H}^n (K^t + \rho B) \, dt.
$$ 
Finally, using the Steiner formula for the convex bodies $K^t$, we obtain
$$
I(f_\rho) = \int_0^{+\infty} \sum_{i = 0}^n \rho^i 
\left({n \atop i} \right) W_i(K^t) \, dt =
\sum_{i=0}^n \rho^i \left({n \atop i} \right) \int_0^{+\infty} W_i(K^t)\,dt,
$$
which is (3.2).
\qed

\subsection{A dual expansion}
One can observe that the functional notion of mean introduced in Definition 3.1 is not linear with respect 
to the sum in $\QC_\alpha$ (unless $\alpha= -\infty$), while this is always the case for the mean width of convex 
bodies. As the latter quantity can be also defined,
up to a dimensional constant, as
$$
\lim_{\rho\to 0^+}\frac{\H^n(B+\rho K)-\H^n(B)}{\rho} \qquad \forall \, K \in \K ^n\, ,
$$
it is natural to ask what happens,
if in place of considering the map $\rho \mapsto I \big (f \oplus \rho \cdot \Theta _\alpha (B)\big)$ as done in the previous section,  
one looks at its ``dual'' map $\rho \mapsto I \big (\Theta _\alpha (B) \oplus \rho \cdot f\big )$. 

Here we focus attention on the case $\alpha = 0$, namely on the class $\QC_0$ of log-concave functions with the
corresponding algebraic operation. As $\Theta _0(B) = \chi _B$,  we set 
\begin{equation}
\Psi(\rho) := I (\chi _B \oplus \rho \cdot f).
\end{equation}
and
$$
\widetilde M(f) := \lim_{\rho\to0^+}
\frac{\Psi (\rho) - \Psi (0)}{\rho} = \Psi'(0^+),
$$ 
whenever the latter limit exists. The first derivative of the mapping 
$\rho \mapsto \Psi(\rho)$ is by construction linear in $f$ 
(exactly as it occurs for the notion of the mean width introduced 
by Klartag and Milman in \cite{Klartag-Milman05}, mentioned in the Introduction). It turns out 
that $\widetilde M(f)$ is finite only when the support of $f$ is compact: 
in this case it can be computed explicitly, and it is given precisely 
by the logarithm of the maximum of $f$ plus the mean width of the 
support of $f$. More precisely we have the following result, which 
is somehow dual to Theorem 3.4. For this reason we call it 
``dual Steiner-type formula''; however we stress that using this 
expression is somehow an abuse, since in this case the function 
$\rho \mapsto \Psi (\rho)$ is {\it not} 
a polynomial in $\rho$. 

\begin{teo}{\rm (Dual Steiner-type formula)}
Let $f \in \QC_0$ and let $\Psi$ be the mapping defined in $(3.7)$. For every $\rho>0$, there holds
\begin{equation}
\Psi (\rho) = \sum_{j=0}^n 
\left({ n \atop j }\right)\rho^{j+1}
\int_0^{+\infty} W_{n-j} (\cl \{f >t\})\, t^{\rho-1} \,dt.
\end{equation}
In particular, setting $K _f:= \{ f >0 \}$, it holds
\begin{equation}
\widetilde M (f) =
\begin{cases}
\kappa_n 
\log(\max_{\sR^n}(f)) + n W_{n-1}(K _f), & \hbox{if} \ K_f \in \K ^n \\
+\infty, & \hbox{otherwise}. 
\end{cases}
\end{equation}
\end{teo}

\vskip5mm
For the proof of Theorem 3.6 the following elementary Lemma is needed.

\begin{lemma}
For every non-increasing function $g:(0,m] \to \R_+$,
$$
\lim_{\rho \to 0^+} \rho \int_0^m g(t) t^{\rho-1} \, dt = 
\lim_{t \to 0^+} g (t) \, .
$$
\end{lemma}

\proof
Set $L: = g(0+)= \lim_{t \to 0^+} g (t)$.  With a change of variable, we have
$$
\rho \int_0^m g(t) t^{\rho-1} \,dt = 
\int_0^{m^\rho} g\big(t^{1/\rho}\big) \,dt.
$$
If $m \geq 1$, write
\begin{equation}
\int_0^{m^\rho} g\big(t^{1/\rho}\big) \,dt =
\int_1^{m^\rho} g\big(t^{1/\rho}\big) \,dt + 
\int_0^1 g\big(t^{1/\rho}\big) \,dt.
\end{equation}
We observe that the first integral in the r.h.s.\ of (3.10) is 
infinitesimal: since $g$ is non-increasing, we have
$$
\int_1^{m^\rho} g\big(t^{1/\rho}\big) \, dt \leq 
g(1) \big(m^\rho - 1\big).
$$
Concerning the second integral in the r.h.s.\ of (3.10), we observe 
that, as $\rho \to 0^+$, the functions $g\big(t^{1/\rho}\big )$ 
do not decrease and converge pointwise to $L$ on $(0,1$). Hence,
$\int_0^1 g\big(t^{1/\rho}\big) \,dt \to + \infty$, by the monotone
convergence theorem. Thus, the statement is proved for $m \geq 1$.

If $0 < m < 1$, for any prescribed $\varepsilon >0$, we have 
$m^\rho > 1- \varepsilon$, for all $\rho$ small enough.
Then regardless of whether $L = + \infty$ or $L < + \infty$, we have
$$
L \geq L m^\rho \geq \int_0^{m^\rho} g\big(t^{1/\rho} \big) \,dt 
 \geq
\int_0^{1 - \varepsilon} g\big(t^{1/\rho}\big) \,dt \to L(1-\varepsilon), \qquad \hbox{ as } \rho \to 0 ^+,
$$
where we used the monotone convergence theorem once more.
The statement then follows by the arbitrariness of $\varepsilon >0$. 
\qed

\bigskip
{\it Proof of Theorem 3.6}. Let us set for brevity 
$f^{(\rho)} := \chi_B\oplus \rho \cdot f$, which in explicit form reads
$$
f^{(\rho)}(z) = 
\sup\big\{f(y)^\rho: x + \rho y = z, \ \|x\| \leq \rho\big\}, 
\qquad \forall \,z \in \R^n.
$$ 
The above definition yields
$$
\big\{f^{(\rho)} > t \big\} = 
\Big\{x: f^\rho \Big(\frac{x}{\rho}\Big) > t\Big\} + B = 
\rho\, \big\{f > t^{1/\rho} \big\} + B.
$$
Therefore,
\begin{equation}
I\big(f^{(\rho)}\big) = 
\int_0^{m^\rho} \H^n\Big(\rho\, \big\{f > t^{1/\rho}\big \} + B\Big)\,dt =
\rho \int_0^m \H^n \Big(\rho\,\big\{f > t\big\} + B\Big) t^{\rho-1}\,dt.
\end{equation}
Letting $\Omega^t = \{f >t\}$ and $K^t = \cl(\Omega^t)$, we have
\begin{equation}
\H^n\Big(\rho \Omega^t + B\Big) =  \H^n\Big(\rho K^t + B\Big) = 
\sum_{j = 0}^n \rho^j \left({n \atop j}\right) W_{n-j} (K^t).
\end{equation}
Inserting (3.12) into (3.11), the equality (3.8) is proved.

Let us now prove (3.9). Set $m = \max_{\sR^n} f$.
We claim that all the terms corresponding to $j \geq 2$ on the 
right-hand side of (3.8) are $o(\rho)$, as $\rho \to 0$. To see this, 
recall that since the functions in $\QC_0$ are log-concave and are 
vanishing at infinity, they must decay exponentially fast (at least).
Hence, there exist constants $\alpha > 0$, $\beta \in \R$, such that
$$
f(x) \leq e^{-(\alpha |x| + \beta)} \qquad \forall x \in \R^n
$$
(see Lemma 2.5 in \cite{CoFr}), which yields
$$
K^t \subseteq \big\{x: e^{-(\alpha |x| + \beta)} \geq t \big\} = 
\Big\{x: |x| \leq - \frac{\beta + \log t}{\alpha}\Big\}.
$$
Letting $R(t) = \max\{0, - \frac{\beta + \log t}{\alpha}\big\}$, we get
$$
\rho^{j+1}
\int_0^m W_{n-j} (\cl \{f>t\})\, t^{\rho - 1} \, dt \leq 
\rho^{j+1} \int_0^m R(t)^j\, t^{\rho - 1} \,dt = 
\rho^j \int_0^{m^\rho} R(t^{1/\rho})^j \ dt \leq C \rho^j,
$$
and the claim is proved.

Next we observe that the terms corresponding to $j=0$ and $j=1$ in 
the sum of (3.8) are given respectively by
$$
\rho \int_0^m W_n(K^t)\, t^{\rho-1} \,dt = 
\kappa_n\,  \rho \int_0^m t^{\rho-1} \,dt = 
\kappa_n m^\rho = \kappa _n m^\rho = I(f_0)\, m^\rho 
$$
(where in the first equality we have exploited the identity  
$W_n(K^t) = \kappa_n$), and by
$$
n \rho^2 \int_0^m W_{n-1}(K^t)\, t^{\rho-1} \,dt. 
$$
Summarizing, we have
$$
I(f^{(\rho)}) = I(f_0)\, m^\rho + 
n \rho^2 \int_0^m W_{n-1}(K^t)\, t^{\rho-1} \,dt + o(\rho),
$$
whence
$$
\frac{I(f^{(\rho)}) - I (f_0)}{\rho} = 
\kappa_n \frac{m^\rho - 1}{\rho} + 
n \rho \int_0^m W_{n-1}(K^t)\, t^{\rho-1} \,dt + \frac{o(\rho)}{\rho}.
$$
In the limit as $\rho \to 0^+$, the first addendum tends to 
$\kappa_n \log m$, whereas the second one tends to $n W_{n-1}(K _f)$ thanks to Lemma 3.7. 
\qed

\section{Generalized Pr\'ekopa-Leindler inequalities}

This section is entirely devoted to the study of generalized versions 
of the Pr\'ekopa-Leindler inequality. More precisely:
in Section 4.1 we prove some variants of such inequality
for functions of one variable; 
in Sections 4.2-4.3 we extend Pr\'ekopa-Leindler's Theorem from the 
usual case of the volume functional to the general case of arbitrary 
monotone concave functionals on $\K^n$ (including as special cases 
the functional quermassintegrals); 
in Section 4.4 we show that this generalized concavity fails 
to be true if one chooses to define the perimeter of quasi-concave 
functions in a different, though apparently natural, way.

\subsection{Variant of Pr\'ekopa-Leindler inequality in dimension one}

Let us return to Theorem 2.6, which we consider here in dimension one
for non-negative functions defined on $(0,+\infty)$. In some 
situations it is desirable to replace the arithmetic mean
$M_1^{(\lambda)}(x,y)$ on the left-hand side of (2.14) by more 
general means $M_\gamma^{(\lambda)}(x,y)$. In the (rather typical) case, 
when $h$ is non-increasing (and if $\gamma < 1$), this would give 
a strengthened one-dimensional variant of this theorem, since 
the hypothesis would be weaker (due to the inequality
$M_\gamma^{(\lambda)}(x,y) \leq M_1^{(\lambda)}(x,y) $).
The case $\gamma = \alpha = 0$ (and hence $\beta = 0$) was 
considered by K. Ball \cite{Ball}, who showed that the hypothesis
\be
h(M_0^{(\lambda)}(x,y)) \geq M_0^{(\lambda)}(f(x),g(y)),
\qquad \forall \, x,y >0,
\en
implies
\be
\int_0^{+\infty} \! h\ \geq \, M_0^{(\lambda)}
\left(\int_0^{+\infty} \!f, \int_0^{+\infty} \!g \right).
\en
Actually, this assertion immediately follows from Prekopa-Leindler's 
Theorem 2.5, when it is applied in one dimension to the functions
$f(e^{-x}) e^{-x}$, $g(e^{-x}) e^{-x}$ and $h(e^{-x}) e^{-x}$.

Below we propose an extension of Ball's observation
to general values $\gamma \leq 1$.

\vskip5mm
\begin{teo} 
Let $\lambda \in (0,1)$,
$\gamma \in [-\infty, 1]$ and $\alpha \in [-\gamma,+\infty]$.
Let $f,g,h$ be non-negative measurable functions on $(0,+\infty)$. If
\be
h(M_\gamma^{(\lambda)}(x,y)) \geq M_\alpha^{(\lambda)}(f(x),g(y)),
\qquad \forall \, x,y >0  \hbox{ such that } f(x)g(y) > 0\,,
\en
then
\be
\int_0^{+\infty} \!h\ \geq \, M_{\beta}^{(\lambda)}
\left(\int_0^{+\infty} \!f, \int_0^{+\infty} \!g \right) \qquad 
\hbox{with} \ \ \beta = \frac{\alpha \gamma}{\alpha + \gamma}\, .
\en
In the extreme cases $\alpha = -\gamma$ and $\alpha = +\infty$, 
the definition of $\beta$ in $(4.4)$ is understood respectively as 
$\beta = -\infty$ and $\beta = \gamma$. In addition, we put 
$\beta = -\infty$ in case $\gamma = - \infty$.
\end{teo}

\vskip5mm
Before giving the proof of Theorem 4.1 let us recall that, 
as a consequence of the generalized H\"older inequality, we have the following elementary inequality: For all $u_1,u_2,v_1,v_2 \geq 0$ and
$\lambda \in (0,1)$, it holds
\be
M_{\alpha_1}^{(\lambda)}(u_1,v_1) M_{\alpha_2}^{(\lambda)}(u_2,v_2)\,
\geq\, M_{\alpha_0}^{(\lambda)}(u_1 u_2,v_1 v_2),
\en
whenever
\be
\alpha_1 + \alpha_2 > 0, \qquad
\frac{1}{\alpha_0} = \frac{1}{\alpha_1} + \frac{1}{\alpha_2}.
\en

\vskip5mm
Inequality (4.5) also holds in the following cases:

\vskip2mm
$\bullet$ \ $\alpha_0 = \alpha_1 = 0$, \ $0 \leq \alpha_2 \leq +\infty$;

$\bullet$ \ $\alpha_0 = \alpha_2 = 0$, \ $0 \leq \alpha_1 \leq +\infty$;

$\bullet$ \ $\alpha_0 = -\infty$, \ $\alpha_1 + \alpha_2 \geq 0$.

\vskip2mm
The latter includes the cases $\alpha_1 = -\infty$, $\alpha_2 = +\infty$
and $\alpha_1 = +\infty$, $\alpha_2 = -\infty$. Clearly, $\alpha_0 > 0$ when $\alpha_1 > 0$ and $\alpha_2 > 0$;
on the other hand if $\alpha_1 < 0 < \alpha_2$ or $\alpha_2 < 0 < \alpha_1$, then
necessarily $\alpha_0 < 0$.

\vskip5mm
{\it Proof of Theorem 4.1.} 

If $\gamma = 1$, we are reduced to Bracamp-Lieb's Theorem 2.6 
in dimension one.

If $\gamma = 0$, then $\beta = 0$ regardless of $\alpha \geq 0$.
But the hypothesis (4.3) is weaker for $\alpha = 0$, and this
case corresponds to Ball's result $(4.1) \Rightarrow (4.2)$.

Hence, we may assume that $-\infty \leq \gamma < 1$, $\gamma \neq 0$. 
Let $-\gamma \leq \alpha \leq +\infty$ with $\gamma > -\infty$.
In terms of the functions
$$
u(x) = f(x^{1/\gamma}), \quad v(x) = g(x^{1/\gamma}), \quad
w(x) = h(x^{1/\gamma})
$$
the hypothesis (4.3) may be rewritten as
\be
w(z) \geq M_\alpha^{(\lambda)}(u(x),v(y)), \qquad 
z = (1 - \lambda)x + \lambda y, \ \ \  \forall \, x,y > 0 \  
\hbox{ such that } \, u (x)v(y) >0\ .
\en
Here and below we omit for brevity the parameter
$\lambda$ and write just $M_\alpha$ instead of $M_\alpha^{(\lambda)}$.

We apply the inequality (4.5) with $\alpha_1 = \alpha$,
$\alpha_2 = \gamma' = \frac{\gamma}{1-\gamma}$, in which case
the condition (4.6) becomes $\alpha + \gamma' > 0$. Using (4.7), it gives
\bee
w(z) z^{1/\gamma'}
 & = &
w(z) \, M_{\gamma'}(x^{1/\gamma'},y^{1/\gamma'}) \\
 & \geq &
M_\alpha(u(x),v(y)) \, M_{\gamma'}(x^{1/\gamma'},y^{1/\gamma'})
 \ \geq \
M_{\alpha_0}(u(x) x^{1/\gamma'},v(y)y^{1/\gamma'}),
\ene
where $\alpha_0$ is defined by
$$
\frac{1}{\alpha_0} = \frac{1}{\alpha} + \frac{1}{\gamma'} =
\frac{1}{\alpha} + \frac{1}{\gamma} - 1.
$$
Here, in case $\alpha = +\infty$, we have $\alpha_0 = \gamma'$,
and in case $\alpha = 0$, one should put $\alpha_0 = 0$
(with constraint $\gamma > 0$ in view of $\alpha + \gamma'>0$).

Thus, the new three functions
$u(x) x^{1/\gamma'}$, $v(x)x^{1/\gamma'}$ and $w(x) x^{1/\gamma'}$
satisfy the condition (2.14) in one-dimensional Brascamp-Lieb's Theorem 
with parameter $\alpha_0$. Hence, if $\alpha_0 \geq -1$, we obtain 
the inequality (2.15) for these functions, that is,
\be
\int_0^{+\infty} w(z) z^{1/\gamma'} dz\ \geq \, M_{\beta}^{(\lambda)}
\left(\int_0^{+\infty} u(x) x^{1/\gamma'} dx,
\int_0^{+\infty} v(y)y^{1/\gamma'} dy\right)
\en
with $\beta = \frac{\alpha_0}{1 + \alpha_0}$. But
$$
\int_0^{+\infty} u(x) x^{1/\gamma'}\,dx =
\int_0^{+\infty} f(x^{1/\gamma})\, x^{1/\gamma - 1}\,dx =
|\gamma|
\int_0^{+\infty} f(x)\,dx,
$$
and similarly for the couples $(v,g)$ and $(w,h)$. In addition,
$$
\beta = \frac{1}{\frac{1}{\alpha_0} + 1} =
\frac{1}{\frac{1}{\alpha} + \frac{1}{\gamma}} =
\frac{\alpha \gamma}{\alpha + \gamma}.
$$
Here, $\beta = \gamma$ for $\alpha = +\infty$, and
$\beta = 0$ for $\alpha = 0$ and $\gamma > 0$, and
$\beta = -\infty$, for $\alpha = -\gamma$.

Thus, (4.8) yields the desired inequality (4.4) of Theorem 4.1,
provided that:

\vskip2mm
$a)$ \ $\alpha + \gamma' > 0$;

$b)$ \ $\alpha_0 \geq -1$.

\vskip2mm
{\it Case} $0 < \gamma < 1$.

Then $\gamma' > 0$. If $\alpha>0$, then $\alpha_0 > 0$, so both
$a)$ and $b)$ are fulfilled. If $\alpha = 0$, then
$\alpha_0 = 0$, so $a)$ and $b)$ are fulfilled, as well. If $\alpha < 0$,
then necessarily $\alpha_0 < 0$ (as already noticed before).
In this case,
$$
\alpha + \gamma' > 0 \Leftrightarrow -\alpha < \gamma' \Leftrightarrow -\frac{1}{\alpha} > \frac{1}{\gamma'} \Leftrightarrow
\frac{1}{\alpha} + \frac{1}{\gamma} < 1.
$$
In addition, since $b)$ may be rewritten as $-\frac{1}{\alpha_0} \geq 1$,
this condition is equivalent to
$-(\frac{1}{\alpha} + \frac{1}{\gamma'}) \geq 1 \Leftrightarrow
\frac{1}{\alpha} + \frac{1}{\gamma} \leq 0 \Leftrightarrow \gamma \geq -\alpha$, which was assumed.

\vskip2mm
{\it Case} $-\infty < \gamma < 0$.

Then $\gamma' < 0$ and $\alpha>0$ to meet $a)$. Again
$\alpha_0 < 0$, so $b)$ may be written as $-\frac{1}{\alpha_0} \geq 1$.
As before, we have
$$
\alpha + \gamma' > 0 \Leftrightarrow \alpha > -\gamma' \Leftrightarrow \frac{1}{\alpha} > -\frac{1}{\gamma'} \Leftrightarrow
\frac{1}{\alpha} + \frac{1}{\gamma} < 1.
$$
In addition, $b)$ is equivalent to
$-(\frac{1}{\alpha} + \frac{1}{\gamma'}) \geq 1 \Leftrightarrow
\frac{1}{\alpha} + \frac{1}{\gamma} \leq 0 \Leftrightarrow \gamma \geq -\alpha$.

\vskip2mm
{\it Case} $\gamma = -\infty$.

This case may be treated by a direct argument. Indeed, necessarily
$\alpha = +\infty$, and the hypothesis (4.3) takes the form
\be
h(\min(x,y)) \geq \max(f(x),g(y)) \qquad \forall \, x, y\;\mbox{such that} \
f(x)g(y) > 0\ .
\en
We may assume that both $f$ and $g$ are not identically zero. Put
$$
a = \sup\{x>0:f(x)>0\}, \qquad b = \sup\{y>0:g(y)>0\}, 
$$
and let for definiteness $a \leq b \leq +\infty$. If $0<x<a$ and 
$f(x) > 0$, one may choose $y \geq x$ such that $g(y) > 0$, and then 
(4.9) gives $h(x) \geq f(x)$. Hence,
\bee
\int_0^{+\infty} h(x)\,dx
 & \geq &
\int_{\{0<x<a, \, f(x)>0\}} h(x)\,dx \\
 & \geq &
\int_{\{0<x<a, \, f(x)>0\}} f(x)\,dx \ = \ \int_0^{+\infty} f(x)\,dx.
\ene
As a result,
$$
\int_0^{+\infty} h(x)\,dx \geq
\min\Big\{\int_0^{+\infty} f(x)\,dx, \int_0^{+\infty} g(x)\,dx\Big\},
$$
which is the desired inequality (4.4) with $\beta = -\infty$.

Theorem 4.1 is now proved. \qed

\subsection{Pr\'ekopa-Leindler inequality for monotone $\gamma$-concave functionals.}

We are now ready to extend Theorem 2.6 by Brascamp and Lieb to general 
monotone $\gamma$-concave set functionals $\Phi$, mentioned in the Introduction. To be more precise,  a functional $\Phi$ defined on the class of all 
Borel subsets of $\R^n$ with values in $[0,+\infty]$ will be said 
to be monotone, if
$$
\Phi(K_0) \leq \Phi (K_1), \quad \hbox{ whenever } \ K_0 \subseteq K_1,
$$
and to be $(\gamma,\lambda)$-concave with parameters 
$\gamma \in [-\infty,+\infty]$ and $\lambda \in (0,1)$, if
\be
\Phi((1 - \lambda) K_0 + \lambda K_1) \geq
M_\gamma^{(\lambda)} \big(\Phi(K_0),\Phi(K_1)\big) \, , 
\en
for all Borel sets $K_0,K_1$ such that $\Phi(K_0) > 0$ and $\Phi(K_1) > 0$.
If (4.10) is fulfilled for an arbitrary 
$\lambda \in (0,1)$, then we simply say that $\Phi$ is $\gamma$-concave.

We always assume that $\Phi(\emptyset) = 0$. In particular, 
the requirement $\Phi(K) > 0$ ensures that $K$ is non-empty.

If $\Phi$ is monotone, we extend it canonically to the class of all 
Borel measurable non-negative functions on $\R^n$ by setting
$$
\Phi(f) = \int_0^{+\infty} \Phi(\{f \geq r\}) \,dr. 
$$

In case $\Phi$ is well-defined only on $\K^n$, the above definition 
remains well-posed in the class of all semi-continuous, quasi-concave
non-negative functions on $\R^n$.

\vskip5mm
\begin{teo}
Let $\Phi$ be a monotone $(\gamma,\lambda)$-concave functional 
on Borel sets of $\R^n$ $($respectively, on $\K^n)$,
with parameters $\gamma \in [-\infty,1]$ and $\lambda \in (0,1)$.
Let $\alpha \in [-\gamma,+\infty]$, and let 
$f,g,h:\R^n \rightarrow [0,+\infty)$
be Borel measurable $($respectively, semi-continuous quasi-concave$)$ 
functions. If 
\be
h((1 - \lambda) x + \lambda y) \geq M_\alpha^{(\lambda)}(f(x),g(y))
\qquad \forall \, x,y \in \R^n \hbox{ such that } f(x)g(y)>0,
\en
then
\be
\Phi(h) \geq M_{\beta}^{(\lambda)}(\Phi(f),\Phi(g)) \quad 
\hbox{ where } \, \beta := \frac{\alpha \gamma}{\alpha + \gamma}\ .
\en
\end{teo}

\vskip5mm
Before giving the proof, several comments on the above statement are 
in order.

\begin{remark}
{\rm
(i) Theorem 2.6 by Brascamp-Lieb can be recast as a special case from 
Theorem 4.2 by taking for the functional $\Phi$ the Lebesgue measure 
on $\R^n$, in which case $\gamma = \frac{1}{n}$.

\vskip2mm
(ii)  In the extreme cases the interpretation of the parameter $\beta$ 
in Theorem 4.2, as well as in the Corollaries hereafter, has to be the 
same as in Theorem 4.1. 

\vskip2mm
(iii) In particular, $\beta = \gamma$ for $\alpha = +\infty$. Thus,  if $f = \chi_{K_0}$, $g = \chi_{K_1}$, and
$h = \chi_{(1 - \lambda) K_0 + \lambda K_1}$, the inequality
(4.11) is fulfilled, and (4.12) gives back the definition of $\gamma$-concavity of $\Phi$. In other words, Theorem 4.2 does 
represent a functional form for the geometric inequality (4.10).

\vskip2mm
(iv) The proof of Theorem 4.2 given below is obtained without using 
an induction argument on the space dimension $n$, but just combining
the $\gamma$-concavity inequality satisfied by assumption by $\Phi$, 
with the one-dimensional functional result stated in Theorem 4.1.

\vskip2mm
(v) If a functional $\Phi$ is monotone and $\gamma$-concave on a given subclass of Borel sets (possibly different than $\K^n$), our proof of 
Theorem 4.2 shows that the implication $(4.11) \Rightarrow (4.12)$
holds true for all Borel measurable functions whose level sets belong 
to the class under consideration.
}
\end{remark}

\bigskip
{\it Proof of Theorem 4.2.} Denote by $K_f(r)$ the level sets
$\{f \geq r\}$, and similarly for $g$ and $h$.
By the hypothesis (4.11), we have the set inclusion
\be
(1-\lambda) K_f(r) + \lambda K_g(s) \subseteq
K_h\big(M_{\alpha}^{(\lambda)}(r,s)\big),
\en
which makes sense and is valid for all $r,s > 0$ such that 
$\Phi(K_f(r))>0$ and $\Phi(K_g(s))>0$. Using (4.13), together 
with the monotonicity and $(\gamma,\lambda)$-concavity assumption 
on $\Phi$, we see that the functions
$$
u(r) := \Phi\big(\{f \geq r\}\big), \quad 
v(r) := \Phi\big(\{g \geq r\}\big), \quad
w(r) := \Phi\big(\{h \geq r\big)\}
$$
satisfy the relation
$$
w\big (M_{\gamma}^{(\lambda)}(r,s)\big ) \geq
M_{\alpha}^{(\lambda)}(u(r),v(s)), \quad \hbox{ whenever } \ u(r)v(s)>0\ .
$$

Therefore, we are in position to apply Theorem 4.1 to the triple
$(u,v,w)$, which yields
$$
\int_0^{+\infty} w(r)\,dr\ \geq \, M_{\beta}^{(\lambda)}
\left(\int_0^{+\infty} u(r)\,dr, \int_0^{+\infty} v(r)\,dr\right)
$$
with $\beta = \frac{\alpha \gamma}{\alpha + \gamma}$.
This is exactly (4.12). \qed

\bigskip
\subsection{Hyperbolic functionals.}

Let us now specialize Theorem 4.2 to an important family of geometric
functionals called hyperbolic or convex. 

\vskip5mm
\begin{definition}{\rm
A monotone functional $\Phi$ defined on the class of all 
Borel subsets of $\R^n$ with values into $[0, + \infty]$ is said to be {\it hyperbolic}, if
\be
\Phi((1 - \lambda) K_0 + \lambda K_1) \geq
\min\big\{\big(\Phi(K_0),\Phi(K_1)\big)\big\}, 
\en
for all $\lambda \in (0,1)$ and for all Borel sets $K_0,K_1$ in $\R^n$ 
such that $\Phi(K_0)>0$ and $\Phi(K_1) > 0$.
}
\end{definition}

\vskip2mm
We adopt a similar definition also if $\Phi$ is defined only on some sublass
of Borel sets, such as $\K^n$.
Thus, hyperbolic functionals are exactly $(-\infty)$-concave
functionals, {\it i.e.}, they satisfy the inequality (4.10) with 
$\gamma = -\infty$. 

\vskip2mm
Apparently, the application of Theorem 4.2 to hyperbolic functionals seems to be not 
so interesting. Indeed, when
$\gamma = -\infty$,  one has $\alpha = +\infty$, in
which case the hypothesis (4.11) considerably restricts the range of
applicability of the resulting inequality (4.12). Nevertheless, the situation is much more favorable if the
hyperbolicity condition (4.14) is combined with some homogeneity property.

\vskip5mm
\begin{definition}{\rm
A functional $\Phi$ defined on the class of all Borel subsets of $\R^n$
(respectively on convex compact sets in $\R^n$) is said to be 
{\it homogeneous of order $\rho$} (with $\rho \in \R\setminus  0$), if
\be
\Phi(\lambda K) = \lambda^\rho\, \Phi(K), 
\en
for all $\lambda > 0$ and for all Borel sets $K$ in $\R^n$ 
(respectively, for all $K \in \K^n$).
}
\end{definition}

\vskip5mm
Combining (4.14) and (4.15) yields the following observation, which is elementary and well-known, especially
for the Lebesgue measure. However, because of its importance, we state it 
separately and in a general setting:

\vskip5mm
\begin{prop}
Any hyperbolic functional $\Phi$, which is homogeneous of order 
$\rho$, is $\gamma$-concave for $\gamma = 1/\rho$. 
\end{prop}

\proof Let $\Phi(K_0) > 0$ and $\Phi(K_1) > 0$.  We have to show that
\be
\Phi(K_0 + K_1) \geq 
\big(\Phi(K_0)^\gamma + \Phi(K_1)^\gamma\big)^{1/\gamma},
\en
If $\Phi(K_0 + K_1) = +\infty$, then (4.16) is immediate. Otherwise, 
$0 < \Phi(K_0) < +\infty$ and $0 < \Phi(K_1) < +\infty$, by the 
monotonicity of $\Phi$. In this case, set
$$
K_0' := \frac{1}{\Phi(K_0)^\gamma}\, K_0 \qquad {\rm and} \qquad 
K_1' := \frac{1}{\Phi(K_1)^\gamma}\, K_1, 
$$
so that, by the homogeneity property (4.15), $\Phi(K_0') = \Phi(K_1') = 1$.
Next, applying the assumption (4.14) to $K_0'$ and $K_1'$, with 
$$
\lambda = \frac{\Phi(K_1)^\gamma}{\Phi(K_0)^\gamma + \Phi(K_1)^\gamma}\, ,
$$
and using once more (4.15), we arrive exactly at the desired 
inequality (4.16).

Finally, being applied to the sets $(1-\lambda)K_0$ and $\lambda K_1$
with arbitrary $\lambda \in (0,1)$, (4.16) turns into (4.10),
expressing the $\gamma$-concavity property of the functional $\Phi$.
\qed
 
\vskip5mm
As a consequence of Proposition 4.6, one may apply Theorem 4.2 to hyperbolic functionals 
$\Phi$, which are homogeneous of order $\rho$, 
as long as $\gamma = \frac{1}{\rho} \leq 1$, that is, when $\rho < 0$ 
or $\rho \geq 1$. In that case, if $\lambda \in (0,1)$, 
$\alpha \in [-\gamma,+\infty]$, and if the functions $f,g,h \geq 0$ 
on $\R^n$ satisfy
\be
h((1 - \lambda) x + \lambda y) \geq M_\alpha^{(\lambda)}(f(x),g(y))
\qquad \forall \, x,y \in \R^n \hbox{ such that } f(x)g(y)>0,
\en
we obtain 
\be
\Phi(h) \, \geq \, M_{\beta}^{(\lambda)}(\Phi(f),\Phi(g)) \quad 
\hbox{ with } \  \beta = 
\frac{\alpha \gamma}{\alpha + \gamma} = \frac{\alpha}{1 + \alpha\rho}.
\en

Similarly as done for Theorem 2.7, one may develop a further generalization
of this statement, involving the means $M_\alpha^{(s,t)}$ for arbitrary $s$ and $t>0$,  not 
necessarily satisfying $s+t = 1$,  and taking in (4.17) the ``optimal''
function
$$
h = s \cdot f \oplus t \cdot g.
$$
Here, the operations $\oplus$ and $\cdot$ are those in ${\mathcal C} _{\alpha}$ for a fixed value $\alpha \geq - \gamma$.
Arguing as before, let for simplicity $\alpha$ be non-zero and finite. 
By its definition, for all $x,y \in \R^n$, the above function $h$ satisfies
$$
h(sx+ty) \geq M_\alpha^{(s,t)}(f(x),g(y)) = (s+t)^{1/\alpha}
\Big(\frac{s}{s+t}\,f(x)^\alpha + \frac{t}{s+t}\,g(y)^\alpha\Big)^{1/\alpha},
$$
which means that the triple $(f,g,\tilde h)$, where
$$
\tilde h(z) := (s+t)^{-1/\alpha}\, h((s+t)\,z),
$$ 
satisfies the hypothesis (4.17) with $\lambda = \frac{t}{s+t}$. 
Hence, we obtain (4.18), {\it i.e.},
\be
\Phi\big(\tilde h\big) \, \geq \, \left[\frac{s}{s+t}\,
\Phi(f)^\beta + \frac{t}{s+t}\,\Phi(g)^\beta\right]^{1/\beta}.
\en
Changing the variable and using the homogeneity property (4.15), we find 
\bee
\Phi\big(\tilde h\big)
 & = &
\int_0^{+\infty} 
\Phi\big(\{z: h((s+t)\,z) \geq (s+t)^{1/\alpha} r\}\big)\,dr \\
 & = &
(s+t)^{-1/\alpha} \int_0^{+\infty} 
\Phi\big(\{z: h((s+t)\,z) \geq r\}\big)\,dr \\
 & = &
(s+t)^{-\rho - 1/\alpha} \int_0^{+\infty} 
\Phi\big(\{z: h(z) \geq r\}\big)\,dr \\
 & = & 
(s+t)^{-\rho - 1/\alpha}\ \Phi(h).
\ene
Taking into account that $\rho = \frac{1}{\alpha} + \frac{1}{\beta}$, the inequality 
(4.19) can be reformulated as in the following statement, where we include the limit 
case $\alpha = +\infty$ as well.

\begin{teo} 
Let $\Phi$ be a hyperbolic functional defined on Borel sets of $\R^n$
$($respectively, on $\K^n)$, which is homogeneous of order $\rho$, with $\rho < 0$ or $\rho \geq 1$. Let  $s, t>0$, let
$\alpha \in [-\frac{1}{\rho},+\infty]$, and let $f, g: \R ^n \to [0, + \infty]$ be 
measurable
$($respectively, semi-continuous quasi-concave$)$ functions. Then
\be
\Phi(s \cdot f \oplus t \cdot g) \geq 
M_\beta^{(s,t)}\big(\Phi(f),\Phi(g)\big), \qquad \hbox{ where }
\beta := \frac{\alpha}{1 + \alpha \rho}\,. 
\en
In case $\alpha = 0$, the restriction
$s+t = 1$ is necessary.  
In the extreme cases
$\alpha = -\frac{1}{\rho}$ and $\alpha = +\infty$, the definition of $\beta$ in $(4.20)$ has to be understood respectively as 
$\beta = -\infty$ and $\beta = \frac{1}{\rho}$.
\end{teo}

\bigskip
Note that the space dimension $n$ is not involved in (4.20).
In particular, when $t=s=1$, such inequality becomes
$$
\Phi(f \oplus g) \, \geq \, 
\bigg[\Phi(f)^{\frac{\alpha}{1+\alpha \rho}} +
\Phi(g)^{\frac{\alpha}{1+\alpha \rho}}\bigg]^{\frac
{1+\alpha \rho}{\alpha}}, \qquad \hbox{ where }
\alpha \neq 0, \ \alpha \geq -\frac{1}{\rho}\, .
$$
In a similar way as already discussed in Section 2.5, this may be viewed as an  
extension to hyperbolic functionals in higher dimensions of the result of Henstock and Macbeath, who 
considered the case of the Lebesgue measure in dimension $n=1$.

\bigskip
As a basic example illustrating Theorem 4.7, we apply it to the
quermassintegrals $\Phi = W_i$, which are known to be
hyperbolic and homogeneous of positive (integer) orders $\rho = n-i$.

\begin{cor} 
Let $i = 0,1,\dots,n-1$. Let  $s, t>0$, let $\alpha \in [-\frac{1}{n-i},+\infty]$,
and let
$f,g$ belong to $\QC_\alpha$. Then 
\be
W_i(s \cdot f \oplus t \cdot g) \geq 
M_\beta^{(s,t)}\big(W_i(f),W_i(g)\big), \qquad 
\beta = \frac{\alpha}{1 + \alpha (n-i)}.
\en
In case $\alpha = 0$, the restriction
$s+t = 1$ is necessary.  
\end{cor}

\bigskip
For $i=1$, we recall that  $n W_1(K)$ represents the perimeter of a set 
$K \in {\cal K}_n$, while according to the co-area formula 
({\it cf.} Remark 3.5), the perimeter of any $C^1$-smooth function 
$f$, vanishing at infinity, can be expressed as the integral
$$
\Per(f) = \int |\nabla f(x)|\,dx.
$$
Hence, in this special case, and for $s+t = 1$, Corollary 4.8 can be rephrased as:

\vskip5mm
\begin{cor}
Let $\lambda \in (0,1)$ and $\alpha \in [-\frac{1}{n-1},+\infty]$, 
$n \geq 2$. Let $f,g,h:\R^n \rightarrow [0,+\infty)$ be 
$C^1$-smooth quasi-concave functions, such that
$h(z) \rightarrow 0$ as $|z| \rightarrow \infty$. If
$$
h((1 - \lambda) x + \lambda y) \geq M_\alpha^{(\lambda)}(f(x),g(y))
\qquad \forall \, x,y \in \R^n \hbox{ such that } f(x) g (y) >0\ ,
$$
then
$$
\int |\nabla h(z)|\,dz \, \geq \, M_{\beta}^{(\lambda)}
\Big(\int |\nabla f(x)|\,dx,\int |\nabla g(y)|\,dy\Big) \quad 
\hbox{ with } \ \beta = \frac{\alpha}{1 + \alpha (n-1)}\ .
$$
\end{cor}

\vskip5mm
Here, the hypothesis that $h$ vanishes at infinity guarantees that 
$f$ and $g$ vanish at infinity, as well. Moreover, the $C^1$-smoothness assumption may be relaxed to the property of being locally Lipschitz.

\subsection{Counterexamples}
Below we show that, choosing a different functional equivalent of the unit ball, 
may lead to a notion of perimeter
which does not satisfy a concavity property like the one stated in 
Corollary 4.9. To be more precise, let us  restrict ourselves to the case $\alpha=0$, namely
to the class $\QC_0$ of log-concave functions, endowed with its corresponding
algebraic structure. Then, for a given function $f\in\QC_0$, 
the definition of the perimeter given in 
Section 2 amounts to
$$
{\rm Per}(f) = \lim_{\rho \to 0^+}
\frac{I(f\oplus\rho\cdot\chi_B)-I(f)}{\rho}\,.
$$
In this definition, one might be willing to replace $\chi_B$ by another log-concave 
function acting as a unitary ball. A natural choice would be the
Gaussian function
$$
e^{-|x|^2/2}\,,
$$
or, more generally,
$$
g_q(x)=e^{-|x|^q/q}\,,
$$
with $q\ge1$. Note that this function tends to $\chi_B$ as $q\to+\infty$. 
In this case one could then define
$$
{\rm Per}_q(f) :=
\lim_{\rho\to0^+}\frac{I(f\oplus\rho\cdot g_q)-I(f)}{\rho},
$$
whenever this limit exists. It was proved in \cite{CoFr} that, under
suitable assumption of smoothness, decay at infinity and strict convexity 
of $f$ (see Theorem 4.5 in \cite{CoFr} for the precise statement), 
the following representation formula holds:
$$
{\rm Per}_q(f)=\frac{1}{p}\int_{\R^n}|\nabla u(x)|^p\, f(x)\,dx\,,
$$
where $p= \frac{q}{q-1}$ is the conjugate H\"older exponent of $q$, 
and $f=e^{-u}$. The aim of this section is to show that 
$
{\rm Per}_q(f)$ 
does not have the same significant properties shown in the previous sections for 
${\rm Per}(f)$, and in
particular it does not verify a generalized Pr\'ekopa-Leindler inequality.

For simplicity, given $f\in\QC_0$ and $p\in(0,\infty)$ let
$$
I_p(f)=\int_{\R^n}|\nabla u(x)|^p\, f(x)\,dx\,.
$$
We want to show that, if $f_0,f_1\in\QC_0$, $t\in[0,1]$, and $f_t:=(1-t)\cdot f_0\oplus t\cdot f_1$,
the following inequality is in general false:
\begin{equation}
I_p(f_t)\ge (I_p(f_0))^{1-t}(I_p(f_1))^t \,. 
\end{equation}
We will consider log-concave functions of the form
\begin{equation}
f(x)=e^{-h_K(x)},
\end{equation}
where $K$ is a convex body in ${\bf R}^n$ and $h_K$ is the support 
function of $K$. We will always assume that $K$ contains the origin 
as {\it interior point}. Since support functions are sub-linear and positively
homogeneous of order one (see \cite{Schneider}), in particular they 
are convex. Thus a function of the form (4.23) is
log-concave (and it is also real-valued and non-negative).

The following result is probably well-known; we include the proof for
the sake of completeness.

\begin{prop}
Let $K_0, K_1$ be convex bodies in ${\bf R}^n$, let $f_0=e^{-h_{K_0}}$,
$f_1=e^{-h_{K_1}}$, let  $t\in [0,1]$, and let $
f_t:=(1-t)\cdot f_0\oplus t\cdot f_1$. 
Then
$$
f_t=e^{-h_{K_0\cap K_1}}\,.
$$
\end{prop}

\proof Set $u_t:=-\log(f_t)$; we want to prove that $u_t=h_{K_0\cap K_1}$. 
For every $z\in{\bf R}^n$, setting for brevity $h_i=h_{K_i}$ for $i=0,1$, 
we have
\begin{eqnarray*}
u_t(z)&=&\inf\{(1-t) h_0(x)+th_1(y)\,:\,(1-t)x+ty=z\}\\
      &=&\inf\{h_0(x)+h_1(y)\,:\,x+y=z\}\,.
\end{eqnarray*}
This means that $u_t$ is the {\it infimal convolution} of $h_0$ and $h_1$.
By Theorem 16.4 in \cite{Ro} we have
$$
(u_t)^*=(h_0)^*+(h_1)^*
$$
where $u^*$ denotes the usual conjugate of convex functions:
$$
u^*(z)=\sup_w [(z,w)-u(w)]\,.
$$
On the other hand, it is easy to verify that, for every convex body $K$,
$$
(h_K)^*=I_K\quad{\rm and}\quad (I_K)^*=h_K
$$
(recall that $I_K$ denotes the indicatrix function of $K$). Hence we have
$$
(u_t)^*=I_{K_0}+I_{K_1}=I_{K_0\cap K_1}=(h_{K_0\cap K_1})^*\,.
$$
The proof is concluded taking the conjugates of the first and the last
function of the above chain of equalities.
\qed

\bigskip
For a convex body $K$ with the origin 
in its interior (and for fixed $p>1$), set
$$
F_p(K) :=I_p(e^{-h_K})= \int_{{\bf R}^n}|\nabla h_K(x)|^p e^{-h_K(x)}\,dx\,.
$$
By Proposition 4.10, inequality (4.22) restricted to functions
of the form (4.23) becomes
\begin{equation}
F_p(K_0\cap K_1)\ge
F_p(K_0)^t\,F_p(K_1)^{1-t}\,,\quad\forall\, K_0, K_1\,,\;\forall\, t\in[0,1]\,.
\end{equation}
The above inequality is in turn equivalent to the the fact that the functional $F_p$
is decreasing with respect to set inclusion, in the class of convex bodies
having the origin as interior point:
\begin{equation}
F_p(K)\ge F_p(K')\,,\quad \forall K\subset K'\,.
\end{equation}
Indeed, taking $K_0=K$, $K_1=K'$ and $t=0$ in (4.24) we get (4.25). 
On the other hand, (4.25) implies that for every $K_0$ and
$K_1$ and for every $t\in[0,1]$,
$$
(F_p(K_0\cap K_1))^t\ge F_p(K_0)^t\,,\quad
(F_p(K_0\cap K_1))^{1-t}\ge F_p(K_1)^{1-t}\,.
$$
Multiplying these inequalities term by term we have (4.24).
In Proposition 4.11 below, we construct examples of convex bodies $K$ and $K'$ 
for which (4.25) is false, under the assumptions $p>1$. 
As an immediate consequence, we obtain that also inequality (4.22) fails to be true for $p>1$.


\begin{prop}
For every $n\ge 2$ and every $p>1$, there exist two convex bodies $K$ 
and $K'$ in $\K ^n$ such that $0\in{\rm int}(K\cap K')$ 
and $F_p(K)< F_p(K')$.
\end{prop}

\begin{cor}
For every $n\ge 2$ and every $p>1$, there exist $f_0,f_1\in\QC_0$ and $t \in [0,1]$ 
such that, if  $f_t: = (1-t) \cdot f_0 \oplus t \cdot f_1$, 
then
\begin{equation*}
\int |\nabla f_t(z)|^p\,dz<
\left(\int |\nabla f_0(z)|^p\,dx\right)^{1-t}
\left(\int |\nabla f_1(z)|^p\,dy\right)^{t}\,.
\end{equation*}
\end{cor}

{\it Proof of Proposition 4.11}. 
We write an arbitrary point $x=(x_1,x_2,\dots,x_n)$ of ${\bf R}^n$ 
in polar coordinates $(r,\theta)=(r,\theta_1,\theta_2,\dots,\theta_n)$:
\begin{equation*}
\left\{
\begin{array}{llllll}
x_1=x_1(r,\theta)=x_1(r,\theta_1,\dots,\theta_{n-1})=r\cos\theta_1\\
x_2=x_2(r,\theta)=x_2(r,\theta_1,\dots,\theta_{n-1})=r\sin\theta_1\cos\theta_2\\
x_3=x_3(r,\theta)=x_3(r,\theta_1,\dots,\theta_{n-1})=r\sin\theta_1\sin\theta_2\cos\theta_3\\
\vdots\\
x_{n-1}=x_{n-1}(r,\theta)=x_{n-1}(r,\theta_1,\dots,\theta_{n-1})=r\sin\theta_1\dots\sin\theta_{n-2}\cos\theta_{n-1}\\
x_n=x_n(r,\theta)=x_n(r,\theta_1,\dots,\theta_{n-1})=r\sin\theta_1\dots\sin\theta_{n-2}\sin\theta_{n-1}\,.
\end{array}
\right.
\end{equation*}
Here
$(r,\theta_1,\dots,\theta_{n-1}) 
\in[0,\infty)\times[0,\pi)^{n-2}\times[0,2\pi)$.
The Jacobian of the mapping $x=x(r,\theta)$ is
$r^{n-1}\sin^{n-2}\theta_1\sin^{n-3}\theta_2\dots\sin\theta_{n-2}$. For
brevity we set $S=[0,\pi)^{n-2}\times[0,2\pi)$. Let us also set
$$
H_K(\theta)=h_K(x(1,\theta))\,,\quad\theta\in S\,.
$$
By the homogeneity of $h_K$ we have
$$
h_K(x(r,\theta))=rH_K(\theta)\,,\quad\forall r\ge0\,,\;\theta\in S\,.
$$
The gradient of $h_K$ is positively homogeneous of order $0$, so that 
$|\nabla h_K(x(r,\theta))|$ does not depend on $r$. Hence we put
\begin{equation}
N_K(\theta) = |\nabla h_K(x(r,\theta))|\,.
\end{equation}
The functional $F_p(K)$ can now be written in the following form
$$
F_p(K)=\int_S N_K(\theta)^p\left(\int_0^\infty r^{n-1}
  e^{-rH_K(\theta)}\,dr\right)\phi(\theta)\,d\theta\,,
$$
where
$$
\phi(\theta)=\sin^{n-2}\theta_1\sin^{n-3}\theta_2\dots\sin\theta_{n-2}\,.
$$
After integration with respect to $r$, we get
\begin{equation}
F_p(K) = (n-1)!
\int_S\frac{N_K(\theta)^p}{H_K(\theta)^n}\phi(\theta)\,d\theta\,.
\end{equation}

Using the above formula, we can immediately deduce counterexamples to
(4.25) for $p>n$. Indeed, from (4.27) we see that $F_p$ is homogeneous
of order $(p-n)$ with respect to homotheties. In particular, if
$\alpha>1$ and $K$ is such that $F_p(K)>0$ (for instance, if $K$ is 
a ball centered at the origin), we have
$$
F_p(\alpha K)=\alpha^{p-n}F_p(K)>F_p(K)\,,
$$
and since $\alpha K\supset K$, this is in conflict with (4.25).

The construction of counterexamples for $p\le n$ is still based on 
(4.27), but it is slightly more involved. We set
$$
K_1=B \qquad \hbox{ and  } \qquad
K_2={\rm conv}(B\cup l\,e_1)\,,
$$
where ${\rm conv}$ denotes the convex hull, $l\ge 1$ and 
$e_1=(1,0,\dots,0)$. We will prove that, for every $p>1$, there is 
a suitable choice of $l$ such that
$F_p(K_1)<F_p(K_2)$.
Since clearly $K_1\subset K_2$,  and the origin is interior to both $K_1$ 
and $K_2$, this will provide a counterexample to (4.25). Note that the body 
$K_2$ is rotationally invariant with respect to the $x_1$-axis, so that 
the function $H_{K_2}$ depends on $\theta_1$ only. With abuse of
notations we write
$$
H_{K_2}(\theta_1,\theta_2,\dots,\theta_{n-1})=H_{K_2}(\theta_1)\,.
$$
More precisely, an explicit expression for $H_{K_2}$ can be written down. 
Let $\phi\in[0,\pi/2]$ be such that
$$
l=\frac{1}{\cos\phi}\,.
$$
Then
\begin{equation*}
H_{K_2}(\theta_1)=
\left\{
\begin{array}{lll}
\displaystyle{\frac{\cos\phi}{\cos\theta_1}}\quad\mbox{if
  $\theta_1\in[0,\phi]$}\,,\\
\\
1\quad\mbox{if  $\theta_1\in[\phi,\pi]$\,.}
\end{array}
\right.
\end{equation*}
Next we have to compute the function $N_{K_2}$. Due to the axial symmetry 
it is not hard to see that the following formula holds
$$
N_{K_2}(\theta) = |\nabla h_{K_2}(x(r,\theta))| =
\sqrt{H_{K_2}^2(\theta_1) +
\left(\frac{dH_{K_2}}{d\theta_1}(\theta_1)\right)^2}\,.
$$
Hence
\begin{equation*}
N_{K_2}(\theta_1)=
\left\{
\begin{array}{lll}
\displaystyle{\frac{1}{\cos\phi}}
\quad\mbox{if
  $\theta_1\in[0,\phi]$}\,,\\
\\
1\quad\mbox{if  $\theta_1\in[\phi,\pi]$\,.}
\end{array}
\right.
\end{equation*}
Now we can compute $F_p(K_2)$. We have
\begin{eqnarray*}
F_p(K_2) & = & (n-1)!
\int_S\frac{N_{K_2}(\theta)^p}{H_{K_2}(\theta)^n}\phi(\theta)\,d\theta\\
      {}\\   
 &=&
2\pi(n-1)!
\int_{[0,\pi)^{n-2}}\frac{N_{k_2}(\theta_1)^p}{H_{K_2}(\theta_1)^n}
\phi(\theta_1,\theta_2,\dots\theta_{n-2})d\theta_1 d\theta_2\dots
d\theta_{n-2}\\
      {}\\ 
 &=& 2\pi(n-1)! \ C(n)
\int_0^\pi\frac{N_{k_2}(\theta_1)^p}{H_{K_2}(\theta_1)^n}\sin^{n-2}(\theta_1)
 \,d\theta_1\,,
\end{eqnarray*}
where
$$
C(n)=\prod_{i=1}^{n-3}\int_0^\pi \sin^i t\,dt\,.
$$

Using the explicit expressions that we have found for $H_{K_2}$ and 
$N_{K_2}$ we obtain
$$
F_p(K_2) = 2\pi(n-1)! \ C(n)
\left[
(\cos\phi)^{n-p}\int_0^\phi\frac{(\sin\theta_1)^{n-2}}{(\cos\theta_1)^n}\,d\theta_1
+
\int_\phi^\pi\sin^{n-2}\theta_1\,d\theta_1
\right]\,.
$$
If $p>1$ the following equality holds
$$
\lim_{\phi\to{\frac{\pi}{2}}^-}(\cos\phi)^{n-p}\int_0^\phi
\frac{(\sin\theta_1)^{n-2}}{(\cos\theta_1)^n}\,d\theta_1 =\infty
$$
and consequently
$$
\lim_{\phi\to{\frac{\pi}{2}}^-}
F_p(K_2)=\infty\,.
$$

Thus, $F_p(K_2)$ can be made arbitrarily large for a suitable choice 
of $\phi$, and in particular, it can be made strictly bigger that 
$F_p(K_1)$ which is independent of $\phi$. \qed

\bigskip

\bigskip
\section{Integral geometric formulas and the valuation property}

In this section we show that the quantities introduced in Definition 3.1
verify integral geometric formulas and a valuation type property, suitably reformulated in the functional case. In both cases, the proofs are straightforward consequences of the definition of the $W_i$'s and the validity of the corresponding properties for convex bodies.

\subsection{Integral geometric formulae}
To begin with, we introduce a notion of {\it projection} for functions, 
which has already been considered in the literature, see for instance
\cite{Klartag-Milman05}. As in Section 2, for $1\le k\le n$ we denote 
by ${\cal L}^n_k$ the set of linear subspaces of $\R^n$ of dimension $k$.
Furthermore, for $L\in{\cal L}^n_k$, we denote by 
$L^\perp\in{\cal L}^n_{n-k}$ the orthogonal complement of $L$ in $\R^n$. 

\begin{definition}
{\rm Let $k\in\{1,\dots,n\}$, $L\in{\cal L}^n_k$ and $f\in\QC$. 
We define the {\it orthogonal projection of $f$ onto $L$} as the function 
$$
f|L\,:\,E\mapsto[0,+\infty]\,,\quad 
f|L\,(x') = \sup\big\{f(x'+y)\,|\,y\in L^\perp\big\}\,.
$$ 
}\end{definition}

\vskip2mm
When $f$ is the characteristic function of a convex body $K \in \K ^n$, 
for any direction $L\in{\mathcal L}^n_k$, the projection $f|L$ agrees 
with the characteristic function of the projection of $K$ onto $H_\xi$.

The following lemma, whose proof follows directly from Definition 5.1, 
shows that the projection of a quasi-concave function is quasi-concave, 
as well. We recall that for $A\subset\R^n$ and $L\in{\mathcal L}^n_k$, 
$A|L$ denotes the orthogonal projection of $A$ onto $L$.

\begin{lemma}
Let $f\in\QC$, $k\in\{1,\dots,n\}$ and $L\in{\mathcal L}^n_k$. 
For every $t\ge0$,
$$
\left\{x'\in L\,:\,f|L(x')>t\right\}=\left\{x\in\R^n\,:\,f(x)>t\right\}|L\,.
$$
\end{lemma}

\vskip2mm
As a consequence of the Cauchy-Kubota formulas for convex bodies, 
Definition 3.1, Lemma 5.2 and Fubini's Theorem, we
have the following result.

\begin{teo}
{\rm (Cauchy-Kubota integral formula for quasi-concave functions)} 
Given $f\in\QC$, for all integers $1 \leq i \leq k \leq n$,
$$
W_i(f)=c(i,k,n)\,\int_{{\mathcal L}^n_k} W_i(f|L_k)\,dL_k\,,
$$
where the constant $c(i,k,n)$ is the same as in in formula $(2.2)$.
\end{teo}

\vskip2mm
As a special case, we consider $i=k=1$, which corresponds to the Cauchy formula. 

\begin{definition} 
{\rm For $\xi \in S^{n-1}$, let $H_\xi$ denote the hyperplane through 
the origin orthogonal to $\xi$. For every $f \in \A$, we define the 
{\it projection of $f$ in the direction $\xi$} as the function defined 
on $H_\xi$ by
$$
(f|\xi)(x') = \sup\big\{f(x' + s\xi):\ s \in \R\big\}, \qquad x' \in H_\xi\ .
$$
}
\end{definition}

\begin{prop}
{\rm (Cauchy integral formula for quasi-concave functions)}
For any $f \in \QC$,
\begin{equation}
\Per(f) = c_n \int_{S^{n-1}} 
\Big\{\int_{H_\xi}(f|\xi)(x')\, d\H^{n-1}(x')\Big\}\,d\H^{n-1}(\xi)\ .
\end{equation}
\end{prop}

\subsection{Valuation property}
The quermassintegrals of convex bodies are known to satisfy the 
following restricted additivity property: For every $i=0,\dots,n$,
\begin{equation}
W_i(K)+W_i(L)=W_i(K\cup L)+W_i(K\cap L)\,,
\end{equation}
for all $K,L\in\K^n$ such that $K\cup L\in\K^n$. A real-valued functional defined on $\K$ for which (5.2) holds is called a valuation. 
The notion of valuation can be transposed into a functional setting,
simply replacing union and intersection by maximum and minimum. 
At this regard, note that if $f$ and $g$ are quasi-concave function, 
then $f\vee g$ is quasi-concave, as well. Here we prove that 
all quermassintegrals of functions in $\QC$ are valuations in the 
above sense.

\begin{prop} 
{\rm (Valuation property)}
Let $f, g \in \QC$ be such that $f \wedge g \in \QC$. Then
for every $i = 0,1,\dots, n-1$,
$$
W_i(f \wedge g) + W_i(f \vee g) = W_i(f) + W_i (g)\ .
$$
\end{prop}

\proof We observe that, for every $t>0$,
$$
\begin{array}{ll}
& \big\{f \wedge g \geq t \big\} = 
\big\{f \geq t \big\} \cap \big\{g \geq t \big\}  \\ \noalign{\smallskip}
& \big\{f \vee g \geq t \big\} = 
\big\{f \geq t \big\} \cup \big\{g \geq t \big \}\ . \end{array}
$$
Since $f,g \in \QC$, one can easily check that also $f \vee g\in \QC$, whereas $f \wedge g \in \QC$, by the assumption. Therefore all the
superlevels appearing in the above equalities belong to $\K^n$, and 
the valuation property (5.2) for the geometric quermassintegrals ensures that
$$
W_i \big(\{f \wedge g \geq t \}\big) + W_i\big(\{f \vee g \geq t\}\big ) 
= W_i \big(\{g \geq t\}\big ) + W_i\big(\{g \geq t\}\big)\ .
$$
Recalling Definition 3.1, the statement follows after integration 
over $(0,+\infty)$. \qed

\bigskip
\section{Functional inequalities}

As we have explicitly defined a notion of the perimeter for 
quasi-concave functions, it is natural to ask for related isoperimetric 
type inequalities. Below, we propose two different kind of inequalities 
in this direction.

\begin{teo} {\rm (Isoperimetric-type inequalities)}
\begin{itemize}
\item [{\rm (i)}] For every $f \in \QC$,
\begin{equation}
\Per(f)\ge n\kappa_n^{1/n}\,\|f\|_{\frac{n}{n-1}}\,.
\end{equation}
\item[{\rm (ii)}] For every $f \in \QC_0$,
\begin{equation}
\Per(f)\ge n I(f) + \mbox{{\rm Ent}}(f)\,,
\end{equation}
\end{itemize}
where
$$
\mbox{{\rm Ent}}(f)=\int f(x)\log f(x)\,dx - I(f)\,\log I(f)\,.
$$
Equality in $(6.1)$ and $(6.2)$ is attained if and only if $f$ is the characteristic function of an arbitrary ball. 
\end{teo}

\vskip2mm
Inequality (6.1) is nothing but the Sobolev inequality in 
$\R^n$ for functions of bounded variation (for which the equality case 
is known to hold iff $f = \chi _B$ up to translations). Actually, 
it holds without the quasi-concavity assumption. 
Inequality (6.2), 
together with the corresponding equality case, can be obtained by 
applying Theorem 5.1 in \cite{CoFr} with $g = \chi _B$.
The isoperimetric inequality (6.1) can naturally be extended to 
other functional quermassintegrals.

\begin{teo} 
For every $f \in \QC$, and for all integers $0 \leq i \leq k \leq n-1$,
\be
W_k(f) \, \geq \, c\, W_i(f^p)^{1/p}, \quad {where} \ \ 
p = \frac{n-i}{n-k}, \ c = \kappa_n^{1-1/p}.
\en
In particular,
\be
W_k(f) \, \geq \, \kappa_n^{k/n}\, \|f\|_{\frac{n}{n-k}}.
\en
Equality in $(6.3)$ and $(6.4)$ is attained if and only if $f$ is the characteristic function of an arbitrary ball. 
\end{teo}

\vskip5mm
Note that inequality (6.4) corresponds to (6.3) in the particular case $i=0$.
Futhermore, taking $k=1$ in (6.4), gives back  the Sobolev
inequality (6.1).

\vskip2mm
\proof The following inequality holds for the quermassintegrals of convex bodies:
$$
W_k(K) \geq c\, W_i(K)^{1/p},
$$
with $c$ and $p$ as in (6.3),
{\it cf.}\ \cite{Schneider}. Applying this bound to the level sets 
$K_f(t) = \{f \geq t\}$ and integrating over $t > 0$, we therefore obtain
\be
W_k(f) \geq c\, \int_0^{+\infty} W_i\big(K_f(t)\big)^{1/p}\,dt.
\en

To further bound from below the integral in (6.5), we use the following
elementary inequality which is commonly applied in the derivation of the Sobolev inequality (6.1), 
see for instance \cite{Burago-Zalgaller}: If $u = u(t)$ is a non-negative,
non-increasing function on $(0,+\infty)$, then for all $p \geq 1$,
$$
\int_0^{+\infty} u(t)^{1/p}\,dt \geq 
\Big(\int_0^{+\infty} u(t)\,dt^p\Big)^{1/p}.
$$
Choosing $u(t) = W_i(K_f(t))$, we see that the $p$-th power of the
integral in (6.5) is greater than or equal to
$$
\int_0^{+\infty} W_i\big(\{f \geq t\}\big)\,dt^p = 
\int_0^{+\infty} W_i\big(\{f^p \geq t^p\}\big)\,dt^p = 
\int_0^{+\infty} W_i\big(\{f \geq t\}\big)\,dt = W_i(f^p).
$$
\qed

\vskip5mm
While we already noticed that the case $k=1$ in (6.4) amounts to the isoperimetric inequality, the case
$k=n-1$ leads to the following functional version of Urysohn's inequality:

\begin{cor} 
For every $f \in \QC$,
\be
M(f) \, \geq \, 2\kappa_n^{-1/n}\,  \|f\|_n.
\en
Equality in $(6.6)$ is attained if and only if $f$ is the characteristic function of an arbitrary ball. 

\end{cor}

\bigskip
For the characteristic functions of convex bodies, (6.6) reduces to
the classical Urysohn's inequality. We point out that, for log-concave functions, a different 
functional version of the Urysohn inequality involving Gaussian 
densities, was earlier proposed by Klartag and Milman 
in \cite{Klartag-Milman05}.
In fact, (6.4) and its particular case (6.6) admits a further 
refinement in terms of radial functions.
Below, for a given $K \in \K ^n$, we denote by $K^*$ the ball 
with the same mean width as $K$.

\begin{teo}
Given $f \in \QC$, denote by $f^*$ the rearrangement of $f$ obtained 
by replacing each of the level sets $\{ f \geq t\}$ by $\{f\geq t\}^*$. 
Then, for every $k = 0,1, \dots, n-1$,
$$
W_k(f) \geq W_k(f^*)\ .
$$
\end{teo}

\proof We have
\bee
W_k(f) & = & \int_0^{+\infty} W_k\big(\{f \geq t\}\big)\,dt \\
 & \geq &
\int_0^{+\infty} W_k\big(\{f \geq t \}\big)^* \, dt \ = \  
\int_0^{+\infty} W_k\big(\{f^* \geq t\}\big)\, dt \ = \ W_k(f^*)\ .
\ene
\qed

\bigskip\bigskip
S.\ Bobkov\par
Department of Mathematics, University of Minnesota\par
306 Vincent Hall 228\par
207 Church Street SE\par
Minneapolis, MN U.S.A. 55455\par\medskip

A.\ Colesanti \par
Dipartimento di Matematica ``U.Dini'', Universit\`a di Firenze \par
Viale Morgagni 67/A \par
50134 Firenze (Italy) \par \medskip

I.\ Fragal\`a \par
Dipartimento di Matematica, Politecnico di Milano \par
Piazza Leonardo da Vinci, 32 \par
20133 Milano (Italy)
\par \medskip

\end{document}